\documentclass[12pt]{amsart}
\usepackage{cite,amsmath}
\newtheorem{theorem}{Theorem}[section]
\newtheorem{lemma}[theorem]{Lemma}
\newtheorem{corollary}[theorem]{Corollary}
\theoremstyle{definition}
\newtheorem{definition}[theorem]{Definition}
\newtheorem{example}[theorem]{Example}
\theoremstyle{remark}
\newtheorem{remark}[theorem]{Remark}

\numberwithin{equation}{section}

\newcommand{\Nat}{{\mathbb N}}	\newcommand{\Int}{{\mathbb Z}}	
\newcommand{\Real}{{\mathbb R}}	\newcommand{\Com}{{\mathbb C}}	
\renewcommand{\Re}{\operatorname{Re}}	\renewcommand{\Im}{\operatorname{Im}}
\newcommand{\res}{\operatorname{res}}	
\newcommand{\ord}{\operatorname{ord}}	

\newcommand{\nulpunt}{{\scriptscriptstyle\circ}}

\newcommand{\String}{{\mathcal L}}	\newcommand{\GS}{\eta}
		\newcommand{\primestring}{{\mathfrak P}}
	
\newcommand{\Dimensions}{{\mathcal D}}	\newcommand{\primitive}{{N}}
\newcommand{\Rerror}{{R}}
\newcommand{\Flow}{\mathcal{F}}
\newcommand{\weight}{\mathfrak{w}}	
\newcommand{\tweight}{\weight_t}
\newcommand{\seq}{\mathfrak{x}}		\newcommand{\orb}{\mathfrak{p}}

\newcommand{\Domega}{\Delta}		\newcommand{\err}{x}

\begin{document}

\title[A Prime Orbit Theorem for Self-Similar Flows]
{A Prime Orbit Theorem for Self-Similar Flows and Diophantine Approximation}

\author{\sc Michel L.~Lapidus} 

\address{Department of Mathematics,
University of California,
Sproul Hall,
Riverside,
California 92521-0135, USA}
\email{lapidus@math.ucr.edu}

\author{\sc Machiel van Frankenhuysen}

\address{Department of Mathematics,
University of California,
Sproul Hall,
Riverside,
California 92521-0135, USA}
\curraddr{Rutgers University,
Department of Mathematics,
110 Frelinghuysen Road,
Piscataway, NJ 08854-8019, USA}
\email{machiel@math.rutgers.edu}

\thanks{Michel L.~Lapidus was partially supported by the National Science Foundation
under grant DMS-9623002 and DMS-0070497 while he was at the MSRI in Spring 2000.}



%\date{}

\begin{abstract}
Assuming some regularity of the dynamical zeta function,
we establish an explicit formula with an error term for the prime orbit counting function of a suspended flow.
We define the subclass of self-similar flows,
for which we give an extensive analysis of the error term in the corresponding prime orbit theorem.
\end{abstract}

\maketitle

\section{Introduction}
\label{S: intro}

In~\cite{wPmP83},
Parry and Pollicott obtain a Prime Orbit Theorem for certain dynamical systems---the so-called `suspension
flows'.
(See also~\cite[Chapter~6]{wPmP90}.)
The first results of this kind were obtained in special cases by Huber~\cite{Huber},
Sinai~\cite{Sinai},
and Margulis~\cite{Margulis},
among others.
See~\cite[2]{wPmP83} and the relevant references therein,
as well as the historical note in~\cite[p.~154]{tBmKcS91}.
Parry and Pollicott derive the first term in the asymptotic expansion of the counting function of prime orbits,
by applying the Wiener-Ikehara Tauberian Theorem to the logarithmic derivative of the dynamical zeta function.
An alternate approach was taken by Lalley in~\cite[2]{sL89},
who considers,
in particular,
the (approximately) self-similar case.
Using a nonlinear extension of the Renewal Theorem,
he shows that in the nonlattice case,
the leading asymptotics are nonoscillatory.
In the lattice case,
the leading asymptotics are periodic,
and it becomes a natural question whether they are constant or nontrivially periodic.

In a recent book~\cite{long},
we have developed a theory of complex dimensions of fractal strings (one-dimensional drums with fractal
boundary,
see~\cite{mLcP93,mLhM95}).
These (geometric) complex dimensions---defined as the poles of the associated geometric zeta function---enable
us to describe the oscillations intrinsic to the geometry or the spectrum of fractal drums,
via suitable `explicit formulas',
obtained in~\cite[Chapter~4]{long}.

In this paper,
we apply these explicit formulas to obtain an asymptotic expansion for the prime orbit counting function of
suspension flows.
The resulting formula involves a sum of oscillatory terms associated with the dynamical complex dimensions of
the flow.
We then focus on the special case of self-similar flows and deduce from our explicit formulas a Prime Orbit
Theorem with error term.
In the lattice case (to be defined below),
the counting function of the prime orbits,
$\psi_\weight (x)$,
has oscillatory leading asymptotics and our explicit formula enables us to give a very precise expression for
this function in terms of multiplicatively periodic functions.
In the nonlattice case (which is the generic case),
the leading term is nonoscillatory and we provide a detailed analysis of the error term.
The precise order of the error term depends on the `dimension free' region of the dynamical zeta function,
as in the classical Prime Number Theorem.
This region in turn depends on properties of Diophantine approximation of the weights of the flow.

For suspension flows,
the dynamical complex dimensions are defined as the poles of the logarithmic derivative of the dynamical zeta
function.
On the other hand,
the geometric complex dimensions of a fractal string are defined in~\cite[2]{mLmvF99} as the poles of the
geometric zeta function,
which coincides with the dynamical zeta function when the string and the flow are self-similar.
Thus the geometric complex dimensions of a self-similar flow only depend on the poles of the corresponding
zeta function,
and they are counted with a multiplicity,
whereas the dynamical complex dimensions of a flow depend on the zeros and the poles of the dynamical zeta
function,
and they usually have no multiplicity.
Due to the fact that the dynamical zeta function of a self-similar flow has no zeros,
the two sets of complex dimensions coincide in this case.

\section{The Zeta Function of a Dynamical System}
\label{S: DS}

Let $N\geq 0$ be an integer and let $\Omega =\{ 1,\dots ,N\}^\Nat$ be the space of sequences over the alphabet
$\{ 1,\dots ,N\}$.
Let $\weight\colon\Omega\rightarrow (0,\infty ]$ be a function,
called the {\em weight.\/}
On $\Omega$,
we have the left shift $\sigma$,
given on a sequence $(a_n)$ by $(\sigma a)_n=a_{n+1}$.
We define the {\em suspended flow\/} $\Flow_\weight$ on the space $[0,\infty )\times\Omega$ as the following
dynamical system (time evolution,
see~\cite[Chapter~6]{wPmP90}):
\begin{gather}
\Flow_\weight (t,a)=\begin{cases}
(t,a)			&\mbox{if }0\leq t<\weight(a),\\
\Flow_\weight (t-\weight(a),\sigma a)	&\mbox{if }t\geq \weight(a).
\end{cases}
\end{gather}
(Note that $\Flow_\weight (t,a)$ may not be defined.
However,
it is always defined on periodic sequences.)
This formalism is seemingly less general than the one introduced in~\cite[Chapter~1]{wPmP90}.
However,
defining $\weight (a)=\infty$ when the sequence $a$ contains a prohibited word of length $2$,
and $e^{-s\infty}=0$,
allows us to deal with the general case.

Given a finite sequence $\seq =a_1,a_2,\dots ,a_l$ of length $l=l(\seq )$,
we let
\begin{gather*}
a=a_1,a_2,\dots ,a_l,a_1,a_2,\dots ,a_l,\dots
\end{gather*}
be the corresponding periodic sequence,
and we define $\sigma\seq =a_2,\dots,a_l,a_1$.
The {\em total weight\/} of the orbit of $\sigma$ on $\seq$ is
\begin{gather}\label{E: total weight}
\tweight (\seq )=\weight (a)+\weight (\sigma a)+\dots +\weight (\sigma^{l-1}a).
\end{gather}
We now define (see~\cite{rB73,dR76} and~\cite[Chapter~5]{wPmP90}):

\begin{definition}\label{D: dynamical zeta}
The {\em dynamical zeta function\/} of $\Flow_\weight$ is defined as
\begin{align}\label{E: dynamical zeta}
\zeta_\weight (s)=\exp\left(\sum_{\seq}\frac{1}{l(\seq )}e^{-s\tweight (\seq )}\right) ,
\end{align}
where the sum extends over all finite sequences $\seq$ of positive length.
\end{definition}

For $N=0$,
the alphabet is empty,
and we interpret $\Flow_\weight$ as the static flow on a point,
and $\zeta_\weight (s)=1$.
Further,
for $N=1$,
we have the dynamical system of a point moving around a circle of length $\tweight (1)=\weight (1,1,\dots )$,
and $\zeta_\weight (s)=(1-e^{-s\tweight (1)})^{-1}$.

We also introduce the logarithmic derivative
\begin{gather}\label{E: Dlog zeta}
-\frac{\zeta_\weight'}{\zeta_\weight}(s)=\sum_\seq\frac{\tweight (\seq )}{l(\seq )}e^{-s\tweight (\seq )}.
\end{gather}
For $N\geq 1$,
this series does not converge for $s=0$.
We assume that~(\ref{E: Dlog zeta}) converges for some value of $s>0$,
and the abcissa of convergence of this series will be denoted by $D$,
the {\em dimension\/} of $\Flow_\weight$.\footnote
{The dimension often coincides with the topological entropy of the flow;
see~\cite[Chapter~5]{wPmP90} and the references therein.}
Clearly,
$D\geq 0$.
Then~(\ref{E: Dlog zeta}) is absolutely convergent for $\Re s>D$.
Moreover,
as in~\cite[2]{mLmvF99},
we assume that there exists a function $S\colon\Real\rightarrow\Real$,
called the {\em screen},
satisfying $S(t)<D$ for every $t\in\Real$,
such that $-\zeta_\weight'/\zeta_\weight$ has a meromorphic extension to a neighborhood of the region
\begin{align}\label{E: window}
W=\{ s=\sigma +it\colon \sigma\geq S(t)\} ,
\end{align}
called the {\em window}.
In Section~\ref{S: PNT},
we will also assume that $-\zeta_\weight'/\zeta_\weight$ satisfies the growth conditions ({\bf H$_1$})
and ({\bf H$_2$}),
to be introduced in Section~\ref{S: explicit}.
We will then say that $\Flow_\weight$ satisfies ({\bf H$_1$}) and ({\bf H$_2$}).

\begin{definition}
The poles of $-\zeta_\weight'/\zeta_\weight (s)$ in $W$ are called the {\em complex dimensions\/} of the flow
$\Flow_\weight$.
The {\em set of complex dimensions\/} of $\Flow_\weight$ in $W$ is denoted by
$\Dimensions_\weight (W)$ or $\Dimensions_\weight$ for short.
\end{definition}

The nonreal complex dimensions of a flow come in complex conjugate pairs~$\omega ,\,\overline{\omega}$
(provided that $W$ is symmetric about the real axis).
If $\zeta_\weight$ has a meromorphic extension to $W$ as well,
then the complex dimensions of $\Flow_\weight$ are simple and they are located at the zeros and poles of
$\zeta_\weight$,
\begin{gather*}
\Dimensions_\weight(W)=\{\omega\in W\colon\zeta_\weight (\omega )=0\text{ or }\infty\} ,
\end{gather*}
and the residue at a complex dimension $\omega$ (i.e.,
$\res (-\zeta_\weight' /\zeta_\weight ;\omega )$) is $-\ord (\zeta_\weight ;\omega )$,
where $\ord (\zeta_\weight ;\omega )=n$ is the order of $\zeta_\weight$ at $\omega$:
$\zeta_\weight (s)=C(s-\omega )^n+O((s-\omega )^{n+1})$.
In general,
the complex dimensions of $\Flow_\weight$ in $W$ are not simple,
and the residues are not necessarily integers.
By abuse of notation,
we write $\ord (\zeta_\weight ;\omega )=\res (\zeta_\weight' /\zeta_\weight ;\omega )$ if
$\zeta_\weight' /\zeta_\weight$ has a meromorphic extension with a simple pole at $\omega$,
even if the residue is not an integer (and consequently,
$\zeta_\weight$ is not analytic at $\omega$).

\subsection{Periodic Orbits,
Euler Product}

A periodic sequence $a$ in $\Omega$ with period $l$,
$a=a_1,\dots ,a_l,a_1,\dots ,a_l,\dots$,
gives rise to the finite orbit
$\{ a,\sigma a,\dots ,\sigma^{l-1}a\}$ of $\sigma$.
It is clear that $l$ is a multiple of the cardinality $\#\{ a,\sigma a,\dots ,\sigma^{l-1}a\}$ of this orbit.

\begin{definition}
A finite sequence $\seq$ is {\em primitive\/} if its length $l(\seq )$ coincides with the length of the
corresponding periodic orbit of $\sigma$.
\end{definition}

We denote by $\sigma\backslash\Omega$ the space of periodic orbits of $\sigma$.
Thus
\begin{gather}
\sigma\backslash\Omega =\left\{\{\sigma^k\seq\colon k\in\Nat\}\colon\seq\mbox{ is a finite sequence}\right\} .
\end{gather}

We reserve the letter $\orb$ for elements of $\sigma\backslash\Omega$.
So $\orb$ will denote a periodic orbit of~$\sigma$,
and we write $\#\orb$ for its length.
The {\em total weight\/} of an orbit $\orb$ is
\begin{gather}
\tweight (\orb )=\sum_{a\in\orb}\weight (a).
\end{gather}

\begin{theorem}[Euler sum]\label{T: Dlog zeta over orbs}
For $\Re s>D,$
we have the following expression for the logarithmic derivative of $\zeta_\weight:$
\begin{gather}\label{E: Dlog zeta over orbs}
-\frac{\zeta_\weight'}{\zeta_\weight}(s)
=\sum_{\orb\in\sigma\backslash\Omega}\sum_{k=1}^\infty\tweight (\orb )e^{-sk\tweight (\orb )},
\end{gather}
where $\orb$ runs through all periodic orbits of\/ $\Flow_\weight .$
\end{theorem}

\begin{proof}
We write the sum in (\ref{E: Dlog zeta}) over the finite sequences $\seq$ as a sum over the primitive
sequences and repetitions of these.
An orbit $\orb$ of $\sigma$ contains $\#\orb$ different primitive sequences of length $\#\orb$,
so we obtain
\begin{align*}
\sum_\seq\frac{\tweight (\seq )}{l(\seq )}e^{-s\tweight(\seq )}
&=\sum_{\seq\text{:primitive}}\sum_{k=1}^\infty\frac{k\tweight (\seq )}{kl(\seq )}e^{-ks\tweight (\seq )}\\
&=\sum_{\orb\in\sigma\backslash\Omega}\#\orb\sum_{k=1}^\infty
\frac{k\tweight (\orb )}{k\#\orb}e^{-ks\tweight (\orb )}.
\end{align*}
The theorem follows.
\end{proof}

\begin{definition}\label{D: psi w}
The following function counts the periodic orbits and their multiples by their total weight:
\begin{gather}\label{E: psi w}
\psi_\weight (x)=\sum_{k\tweight (\orb )\leq\log x}\tweight (\orb ).
\end{gather}
\end{definition}

The function $\psi_\weight(x)$ is the counterpart of $\psi (x)=\sum_{p^k\leq x}\log p$,
which counts prime powers $p^k$ with a weight $\log p$;
see Example~\ref{E: Riemann-von Mangoldt}.

\begin{corollary}\label{C: logarithmic derivative and psi}
We have the following relation between $\zeta_\weight'/\zeta_\weight$ and $\psi_\weight :$
\begin{gather}\label{E: logarithmic derivative and psi}
-\frac{\zeta_\weight'}{\zeta_\weight}(s)=\int_0^\infty x^{-s}d\psi_\weight (x),
\end{gather}
for $\Re s>D$.
\end{corollary}

The integral on the right-hand side of~(\ref{E: logarithmic derivative and psi}) is a Riemann-Stieltjes
integral associated with the monotonic function $\psi_\weight$.

\begin{corollary}[Euler product]\label{T: suspended flow Euler product}
The function $\zeta_\weight(s)$ has the following expansion as a product\/{\rm :}
\begin{gather}\label{E: suspended flow Euler product}
\zeta_\weight(s)=\prod_{\orb\in\sigma\backslash\Omega}\frac{1}{1-e^{-s\tweight (\orb )}},
\end{gather}
where $\orb$ runs over all periodic orbits of\/ $\Flow_\weight .$
The product converges for $\Re s>D$.
\end{corollary}

\begin{proof}
In (\ref{E: Dlog zeta over orbs}),
we sum over $k$ to obtain
\begin{gather}
\frac{\zeta_\weight'}{\zeta_\weight}(s)
=-\sum_{\orb\in\sigma\backslash\Omega}
\frac{\tweight (\orb )e^{-s\tweight (\orb )}}{1-e^{-s\tweight (\orb )}}
=-\sum_{\orb\in\sigma\backslash\Omega}\frac{d}{ds}\log\left( 1-e^{-s\tweight (\orb )}\right) .
\end{gather}
The theorem then follows upon integrating and taking exponentials.
\end{proof}

In Section~\ref{S: PNT},
we combine the above Euler product representation of $-\zeta_\weight'/\zeta_\weight$ with our explicit
formulas of Section~\ref{S: explicit} to derive a Prime Orbit Theorem for primitive periodic orbits.

\begin{remark}\label{R: psi instead of pi}
We use $\psi_\weight$ instead of the more direct counting function
$$
\pi_\weight (x)=\sum_{\tweight (\orb )\leq\log x}1.
$$
However,
setting $\theta_\weight (x)=\sum_{\tweight (\orb )\leq\log x}\tweight (\orb )$,
so that
$\psi_\weight(x)=\theta_\weight(x)+\theta_\weight(x^{1/2})+\theta_\weight(x^{1/3})
+\dots$ and
$\theta_\weight ({x})=\psi_\weight (x)+O\left(\sqrt{\psi_\weight (x)}\right)$,
as $x\to\infty$,
we find
$$
\pi_\weight (x)=\int_0^x\frac{1}{\log t}\, d\theta_\weight (t)
=\frac{\theta_\weight (x)}{\log x}+\int_0^x\frac{\theta_\weight (t)}{\log^2 t}\,\frac{dt}{t},
$$
from which it is easy to derive the corresponding theorems for $\pi_\weight$ from those for~$\psi_\weight$.
\end{remark}

\subsection{Self-Similar Flows}
\label{S: ssd}

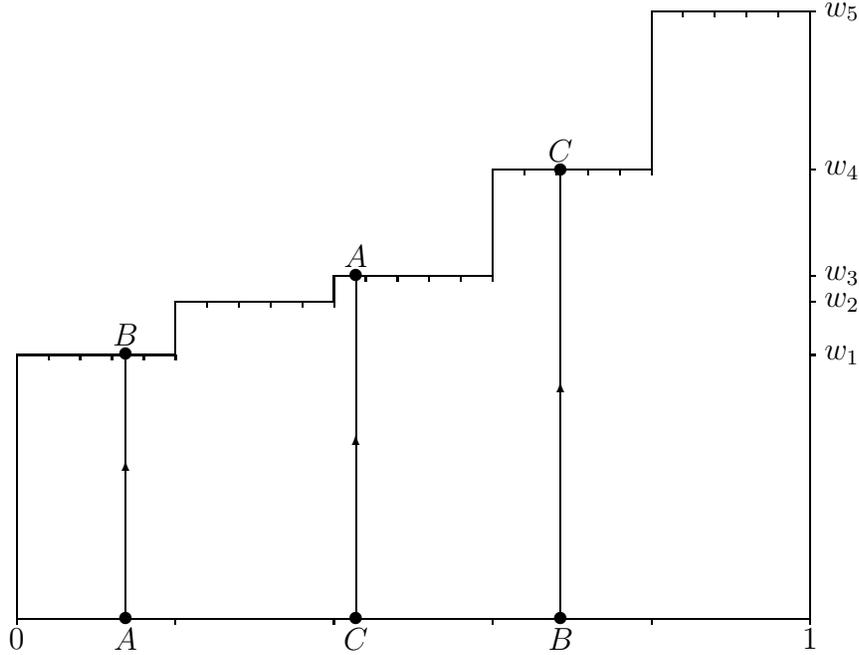
\begin{figure}[tb]
\begin{picture}(324,250)(-10,-10)
%\unitlength .9pt
\put(0,0){\line(1,0){300}}\put(0,0){\line(0,1){100}}\put(300,0){\line(0,1){230}}
\multiput(0,0)(60,0){6}{\line(0,-1){2}}
\put(0,-7.5){\makebox(0,0){$0$}}
\put(300,-7.5){\makebox(0,0){$1$}}

\put(300,100){\put(2,0){\line(-1,0){2}}\put(5.5,0){\makebox(0,0)[l]{$w_1$}}
\put(-300,0){\line(1,0){60}\multiput(0,0)(-12,0){6}{\line(0,-1){2}}}}
\put(300,120){\put(2,0){\line(-1,0){2}}\put(5.5,0){\makebox(0,0)[l]{$w_2$}}
\put(-240,0){\line(1,0){60}\multiput(0,0)(-12,0){6}{\line(0,-1){2}}}}
\put(300,130){\put(2,0){\line(-1,0){2}}\put(5.5,0){\makebox(0,0)[l]{$w_3$}}
\put(-180,0){\line(1,0){60}\multiput(0,0)(-12,0){6}{\line(0,-1){2}}}}
\put(300,170){\put(2,0){\line(-1,0){2}}\put(5.5,0){\makebox(0,0)[l]{$w_4$}}
\put(-120,0){\line(1,0){60}\multiput(0,0)(-12,0){6}{\line(0,-1){2}}}}
\put(300,230){\put(2,0){\line(-1,0){2}}\put(5.5,0){\makebox(0,0)[l]{$w_5$}}
\put(-60,0){\line(1,0){60}\multiput(0,0)(-12,0){6}{\line(0,-1){2}}}}

\put(60,100){\line(0,1){20}}
\put(120,120){\line(0,1){10}}
\put(180,130){\line(0,1){40}}
\put(240,170){\line(0,1){60}}

\put(41.13,0){\put(0,-7.5){\makebox(0,0){$A$}}\put(0,0){\makebox(0,0){$\bullet$}}
\line(0,1){100}\vector(0,1){60}\put(0,107.5){\makebox(0,0){$B$}}\put(0,100){\makebox(0,0){$\bullet$}}}
\put(205.65,0){\put(0,-7.5){\makebox(0,0){$B$}}\put(0,0){\makebox(0,0){$\bullet$}}
\line(0,1){170}\vector(0,1){90}\put(0,177.5){\makebox(0,0){$C$}}\put(0,170){\makebox(0,0){$\bullet$}}}
\put(128.23,0){\put(0,-7.5){\makebox(0,0){$C$}}\put(0,0){\makebox(0,0){$\bullet$}}
\line(0,1){130}\vector(0,1){70}\put(0,137.5){\makebox(0,0){$A$}}\put(0,130){\makebox(0,0){$\bullet$}}}
\end{picture}
\caption{A self-similar flow,
$N=5$,
with the orbit of $17/124$.}
\label{F: ss flow}
\end{figure}

A self-similar flow is best viewed as the following dynamics on the region of Figure~\ref{F: ss flow}.
A point $x=x_1N^{-1}+x_2N^{-2}+\dots =.x_1x_2\dots$ on the unit interval moves vertically upward with unit
speed until it reaches the graph,
at which moment it jumps to $\{ Nx\} =Nx-[Nx]=.x_2x_3\dots$,
the fractional part of $Nx$,
and continues from there.
In Figure~\ref{F: ss flow},
$N=5$,
and the expansions of $A,B$ and $C$ in base $5$ are $A=17/124=.\overline{032}$,
$B=85/124=.\overline{320}$,
$C=53/124=.\overline{203}$.

\begin{definition}
A flow $\Flow_\weight$ is {\em self-similar\/} if\/ $N\geq 2$ and the weight function\/ $\weight$ depends only
on the first letter of the sequence on which it is evaluated.
We then put
\begin{align}\label{E: wj}
w_j=\weight(j,j,j,\dots ),
\end{align}
and
\begin{align}\label{E: rj}
r(\seq )=e^{-\tweight(\seq )},\qquad r_j=e^{-w_j}=r(j,j,j,\dots ),
\end{align}
for $j=1,\dots ,N$.
The numbers $r_j$ are called the {\em scaling ratios\/} of $\Flow_\weight$.
\end{definition}

Note that $0<r_j<1$.
We will assume that the {\em weights\/} $w_j=\log r_j^{-1}$ are ordered in increasing order,
$0<w_1\leq w_2\leq\dots\leq w_N$,
so that $1>r_1\geq r_2\geq\dots\geq r_N>0$.
When $N=2$,
the flow is called a {\em Bernoulli flow}.
Such flows play an important role in ergodic theory (see~\cite[Chapters 2, 6 and 8]{tBmKcS91}).

\begin{theorem}\label{T: dynamical}
The dynamical zeta function associated with a self-similar flow has a meromorphic continuation to the whole
complex plane,
given by
\begin{gather}\label{E: dynamical}
\zeta_\weight(s)=\frac{1}{1-\sum_{j=1}^N r_j^s}.
\end{gather}
Its logarithmic derivative is given by
\begin{gather}\label{E: Dlog dynamical}
-\frac{\zeta_\weight'}{\zeta_\weight}(s)=\frac{\sum_{j=1}^N w_jr_j^s}{1-\sum_{j=1}^N r_j^s}.
\end{gather}
The dimension\/ $D>0$ of the flow is the unique real solution of the equation\/ $1=\sum_{j=1}^N r_j^s$.
\end{theorem}

\begin{proof}
The sum over periodic sequences of fixed length $l$ can be computed as follows:
\begin{align*}
\sum_{\seq :\, l(\seq )=l}r(\seq )^s
&=\sum_{a_1=1}^N\sum_{a_2=1}^N\dots\sum_{a_l=1}^N r_{a_1}^s\dots r_{a_l}^s\\
&=\left( r_1^s+\dots +r_N^s\right)^l.
\end{align*}
Hence,
for $\Re s>D$,
the sum over all periodic sequences is equal to
\begin{gather*}
\sum_{l=1}^\infty\frac{1}{l}
\sum_{\seq :\, l(\seq )=l}r(\seq )^s=\sum_{l=1}^\infty\frac{1}{l}\left( r_1^s+\dots +r_N^s\right)^l
=-\log\biggl( 1-\sum_{j=1}^N r_j^s\biggr) .
\end{gather*}
The theorem follows upon exponentiation and analytic continuation.
Since the function $1-\sum_{j=1}^N r_j^s$ is holomorphic,
$\zeta_\weight$ is meromorphic.
\end{proof}

\begin{remark}\label{R: window is C}
Because of Theorem~\ref{T: dynamical},
for a self-similar flow we can take the full complex plane for the window,
$W=\Com$;
in that case,
there is no screen.
However,
in applying our explicit formulas,
we sometimes choose a screen to obtain information about the error of an approximation.
\end{remark}

\begin{corollary}\label{C: cxd equation}
The set of complex dimensions $\Dimensions_\weight =\Dimensions_\weight (\Com )$ of the self-similar flow
$\Flow_\weight$ is the set of solutions of the equation
\begin{align}\label{E: cxd equation}
\sum_{j=1}^N r_j^\omega =1,
\qquad\omega\in\Com .
\end{align}
Moreover,
the complex dimensions are simple {\rm (}that is,
the pole of\/ $-\zeta_\weight'/\zeta_\weight$ at $\omega$ is simple\/{\rm ).}
The residue at\/ $\omega$ equals $-\ord (\zeta_\weight ;\omega ).$
\end{corollary}

\subsubsection{Connection with Self-Similar Fractal Sets}

Given an open interval $I$ of length~$L$,
we construct a self-similar one-dimensional fractal set~$\String$ with scaling ratios
$r_1,\dots , r_N$.
Such a set is called a {\em fractal string\/} (see~\cite[2]{mLcP93,mLhM95,mLmvF99}).
The following construction is reminiscent of the construction of the Cantor set.
Let $N$ scaling factors $r_1, r_2,\dots , r_N$ be given ($N\geq 2$),
with
$$
1>r_1\geq r_2\geq\ldots\geq r_N>0.
$$
Assume that
\begin{align}\label{E: ss R}
R:=\sum_{j=1}^N r_j<1.
\end{align}
Subdivide $I$ into intervals of length $r_1L,\dots ,r_NL$.
The remaining piece of length $(1-R)L$ is the first member of the string,
denoted by $l_1$,
also called the first length in Remark~\ref{R: special values} below.
Repeat this process with the remaining intervals,
to obtain $N$ new lengths $l_2,\dots ,l_{N+1}$ in the next step,
and $N^{k-1}$ new lengths in the $k$-th step.
As a result,
we obtain a self-similar string $\String$ consisting of intervals of length
$
L(1-R)r_1^{k_1}\ldots r_N^{k_N}$
($k_1,\dots ,k_N\in\Nat$),
and a sequence $l_1\geq l_2\geq l_3\geq\dots$ of positive numbers,
called the {\em lengths\/} of the string.
We let $\zeta_\String (s)=\sum_{j=1}^\infty l_j^s$,
the {\em geometric zeta function\/} of $\String$ (see~\cite[2]{mLmvF99}).

\begin{theorem}\label{T: ss analytic continuation}
Let\/ $\String$ be a self-similar string,
constructed as above with scaling ratios $r_1=e^{-w_1},\dots ,r_N=e^{-w_N}.$
Then the geometric zeta function
of this string has a meromorphic continuation to the whole complex plane,
given by
\begin{equation}\label{E: geoz ss}
\zeta_\String (s)=(L(1-R))^s\zeta_\weight(s),\quad\mbox{for }s\in\Com .
\end{equation}
Here,
$L$ is the total length of\/ $\String ,$
and\/ $R$ is given by~{\rm (\ref{E: ss R}).}
\end{theorem}

This follows from Theorem~\ref{T: dynamical} combined with~\cite[Theorem~2.3, p.~25]{long}.

\begin{remark}\label{R: special values}
For a self-similar string,
the total length of $\String$ is also the length of the initial interval $I$ in the above construction.
We can always normalize a self-similar string in such a way that $\zeta_\String =\zeta_\weight$
(equivalently,
that the first length of $\String$ is~$1$),
by choosing
$
L(1-R)=1.
$
This does not affect the complex dimensions of the string.

Note that we need to assume that $R=\sum_{j=1}^Nr_j<1$,
which corresponds to a lower bound on the weights $w_j=-\log r_j$.
There is no analogue of this condition for general suspended flows.
\end{remark}

\begin{remark}
The Euler product does not seem to have a clear geometric interpretation in the language of fractal strings.
There is,
however,
a natural self-similar flow on $\String$:
the flow
\begin{gather}
\Flow_\String (t,j,x)=\begin{cases}
(0,j,xe^t)	&\mbox{if }xe^t<l_j,\\
\Flow_\String (t-\log l_j,j,1)&\mbox{otherwise}.
\end{cases}
\end{gather}
The lengths $l_j$ correspond to the periodic sequences $\seq$ of the flow $\Flow_\weight$ via the formula
\begin{align}
l_j=\prod_{k=0}^{l(\seq )-1}r(\sigma^k\seq ).
\end{align}
\end{remark}

\begin{remark}[Geometric and dynamical complex dimensions]
In~\cite{long},
the geometric complex dimensions of a fractal string are defined as the poles of its geometric zeta function.
Thus the complex dimensions are counted with a multiplicity,
and the zeros of the geometric zeta function are unimportant.
On the other hand,
the dynamical complex dimensions are defined as the poles of the logarithmic derivative of the dynamical zeta
function.
Thus the complex dimensions are simple,
and both the zeros and the poles of the dynamical zeta function are counted.
For self-similar flows,
the dynamical zeta function and the geometric zeta function of the corresponding string coincide (up to
normalization),
and this zeta function has no zeros.
Hence,
as sets (without multiplicity),
the geometric and dynamical complex dimensions coincide for self-similar flows and strings.
\end{remark}

\begin{remark}[Higher-dimensional case]
We have discussed above the case of fractal strings (i.e.,
the one-dimensional case) because it is the one studied in most detail in~\cite{long}.
However,
it is clear that our results can be applied to higher-dimensional self-similar fractals~\cite{kF90,Mandelbrot}
as well.
This allows us to obtain information about the symbolic dynamics of self-similar fractals.
On the other hand,
as in the previous remark,
it does not give information about the actual geometry of such fractals.
\end{remark}

\subsection{The Lattice and Nonlattice Case}
\label{S: l nl}

Let $\Flow_\weight$ be a self-similar flow.
Recall that $\weight$ depends only on the first symbol and $w_j=\weight (j,j,\dots )$ for $j=1,\dots ,N$.
Consider the subgroup $G$ of $\Real$ generated by these weights,
$G=\sum_{j=1}^N\Int w_j$.

\begin{definition}\label{D: l nl}
The case when $G$ is dense in $\Real$ is called the {\em nonlattice case}.
We then say that $\Flow_\weight$ is a {\em nonlattice flow}.

The case when $G$ is not dense (and hence discrete) in $\Real$ is called the {\em lattice case}.
We then say that $\Flow_\weight$ is a {\em lattice flow}.
In this situation there exists a unique positive real real number $w$,
called the {\em generator\/} of the flow,
and positive integers $k_1,\dots ,k_N$ without common divisor,
such that
$
1\leq k_1\leq\dots\leq k_N
$
and
\begin{align}
w_j={k_j}w,
\end{align}
for $j=1,\dots ,N$.
\end{definition}

The generator of $\Flow_\weight$ generates the flow in the sense that the weight of every periodic sequence of
$\Flow_\weight$ is an integer multiple of $w$.
\medskip

We introduce a real number $D_0$ as follows:
Let~$m$ be the number of integers $j$ in $1,\dots ,N$ such that\/ $r_j=r_N,$
and let $D_0\in\Real$ be defined by
\begin{gather}\label{E: D0}
1+\sum_{j=1}^{N-m}r_j^{D_0}=mr_N^{D_0}.
\end{gather}
The dynamical complex dimensions of a self-similar flow are described in the following theorem.
For brevity,
we will usually refer to them as the complex dimensions of $\Flow_\weight$.

\begin{theorem}\label{T: struc cxd}
Let\/ $\Flow_\weight$ be a self-similar flow of dimension $D$ and with scaling ratios
$1>r_1\geq\dots\geq r_N>0.$
Then the value $s=D$ is the only complex dimension of\/ $\Flow_\weight$ on the real line,
all complex dimensions are simple,
and the residue at a complex dimension {\rm (}i.e.,
$\res (-\zeta_\weight' /\zeta_\weight ;\omega ))$ is a positive integer.
The set of complex dimensions in $\Com$ {\rm (}see Remark\/~{\rm\ref{R: window is C})} of\/ $\Flow_\weight$ is
contained in the bounded strip\/ $D_0\leq\Re s\leq D\colon$
\begin{align}
\Dimensions_\weight =\Dimensions_\weight (\Com )\subset\left\{ s\in\Com\colon D_0\leq\Re s\leq D\right\} .
\end{align}
It is symmetric with respect to the real axis and infinite,
with density bounded by
\begin{align}\label{E: density}
\#\left(\Dimensions_\weight\cap\{\omega\in\Com\colon |\Im\omega |\leq T\}\right)\leq\frac{w_N}{\pi}T+O(1) ,
\end{align}
as $T\to\infty$.

In the {\em lattice case,\/}
$\zeta_\weight(s)$ is a rational function of\/ $e^{-ws},$
where $w$ is the generator of\/ $\Flow_\weight.$
So,
as a function of\/ $s,$
it is periodic with period\/ $2\pi i/w.$
The complex dimensions\/ $\omega$ are obtained by finding the complex solutions\/ $z$ of the polynomial
equation\/
{\rm (}of degree~$k_N$\/{\rm )}
\begin{align}\label{E: algebraic}
\sum_{j=1}^N z^{k_j}=1,\quad\mbox{with } e^{-w\omega}=z.
\end{align}
Hence there exist finitely many poles $\omega_1(=D),\omega_2,\dots ,\omega_q,$
such that
\begin{align}\label{E: lines of cxd}
\Dimensions_\weight =\{\omega_u+2\pi in/w\colon n\in\Int , u=1,\dots ,q\} .
\end{align}
In other words,
the poles lie periodically on finitely many vertical lines,
and on each line they are separated by $2\pi /w$.
The residue of the complex dimensions corresponding to one value of\/ $z=e^{-w\omega}$ is the
multiplicity of\/ $z$ as a solution of\/ {\rm (\ref{E: algebraic}).}

In the {\em nonlattice case,\/}
$D$ is simple and is the unique pole of\/ $\zeta_\weight$ on the line $\Re s=D.$
Further,
there is an infinite sequence of complex dimensions of\/~$\Flow_\weight$ coming arbitrarily close\/
{\rm (}from the left\/{\rm )} to the line $\Re s=D.$
There exists a screen $S$ to the left of the line $\Re s=D,$
such that\/ $-\zeta_\weight'/\zeta_\weight$ satisfies {\rm (\bf H$_1$\rm )} and\/ {\rm ({\bf H$_2$})} with\/
$\kappa =0$\/ {\rm (}see Equations\/~{\rm (\ref{E: H1})} and\/~{\rm (\ref{E: H2})} below\/{\rm ),}
and the residue of\/ $-\zeta_\weight'/\zeta_\weight$ at the pole $\omega$ in $W$ is equal to $1.$
Finally,
the complex dimensions of\/ $\Flow_\weight$ can be approximated\/
{\rm (}via an explicit algorithm,
as described in\/~{\rm\cite[\S2.6]{long})}
by the complex dimensions of a sequence of lattice strings,
with smaller and smaller generator.
Hence the complex dimensions of a nonlattice string have an almost periodic structure.
\end{theorem}

\begin{corollary}\label{C: ss has infty cxd}
Every self-similar flow has infinitely many complex dimensions with positive real part.
\end{corollary}

\begin{proof}[Proof of Theorem~\ref{T: struc cxd}]
For a proof of these facts,
see~\cite[Theorem~2.13, pp.~37--40]{long}.
The density estimate~(\ref{E: density}) follows from the fact that the right-hand
side of~(\ref{E: density}) gives the asymptotic density of the number of poles of
$\zeta_\weight$,
counted {\em with\/} multiplicity.
The $O(1)$ estimate improves~\cite[Theorem~2.22, p.~47]{long}.
It is proved in~\cite{mLmvF01}.
\end{proof}

\subsection{Examples of Complex Dimensions of Self-Similar Flows}

\begin{example}[The Cantor flow]\label{E: Cantor Flow}
This is the self-similar flow on the alphabet $\{ 0,1\}$,
with two equal weights $w_1=w_2=\log 3$.
It has $2^n$ periodic sequences of weight $n\log 3$,
for $n=1,2,\dots$.
The dynamical zeta function of this flow is
\begin{align}
\zeta_{\rm CF}(s)=\frac{1}{1-2\cdot 3^{-s}}.
\end{align}
After taking the logarithmic derivative,
one finds that the dynamical complex dimensions are the solutions of the equation
\begin{align}\label{E: algebraic Cantor}
2\cdot 3^{-\omega}=1\qquad (\omega\in\Com ).
\end{align}
We find
\begin{align}
\Dimensions_{\rm CF}=\left\{ D+\frac{2\pi i}{w}n\colon n\in\Int\right\} ,
\end{align}
with $D=\log_3 2$ and $w=\log 3$.
\end{example}

\begin{example}[The Fibonacci flow]\label{E: Fibonacci Flow}
Next we consider a self-similar flow with two lines of complex dimensions.
The {\em Fibonacci flow\/} is the flow Fib on the alphabet $\{ 0,1\}$ with weights $w_1=\log 2$,
$w_2=2\log 2$.
Its periodic sequences have weight $\log 2, 2\log 2,\dots ,n\log 2, \dots$,
with multiplicity respectively $1,2,\dots ,F_{n+1},\dots ,$
the Fibonacci numbers.
Recall that these numbers are defined by the following recursive equation:
\begin{align}\label{E: Fibonacci numbers}
F_{n+1}=F_n+F_{n-1},\mbox{ with }F_0=0,\ F_1=1.
\end{align}

The dynamical zeta function of the Fibonacci flow is
\begin{align}
\zeta_{\rm Fib}(s)=\frac{1}{1-2^{-s}-4^{-s}}.
\end{align}
The complex dimensions are found by solving the quadratic equation
\begin{align}\label{E: algebraic Fibonacci}
(2^{-\omega})^2+2^{-\omega}=1\qquad (\omega\in\Com ).
\end{align}
We find $2^{-\omega}=\left( -1+\sqrt{5}\right) /2=\phi^{-1}$ and $2^{-\omega}=-\phi$,
where $\phi =(1+\sqrt{5})/2$ is the golden ratio.
Hence
\begin{align}
\Dimensions_{\rm Fib}=\left\{ D+\frac{2\pi i}{w}n\colon n\in\Int\right\}
\cup\left\{ -D+\frac{2\pi i}{w}(n+1/2)\colon n\in\Int\right\} ,
\end{align}
with $D=\log_2\phi$ and $w=\log 2$.
\end{example}

\begin{example}[The Golden flow]\label{E: golden flow}
We consider the nonlattice flow GF with weights $w_1=\log 2$ and $w_2=\phi\log 2$,
where $\phi =(1+\sqrt{5})/2$ is the golden ratio.
We call this flow the {\em golden flow\/}.
Its dynamical zeta function is
\begin{align}\label{E: geo golden}
\zeta_{\rm GF}(s)=\frac{1}{1-2^{-s}-2^{-\phi s}},
\end{align}
and its complex dimensions are the
solutions of the transcendental equation
\begin{align}\label{E: cxd golden}
2^{-\omega}+2^{-\phi\omega}=1\qquad (\omega\in\Com ).
\end{align}

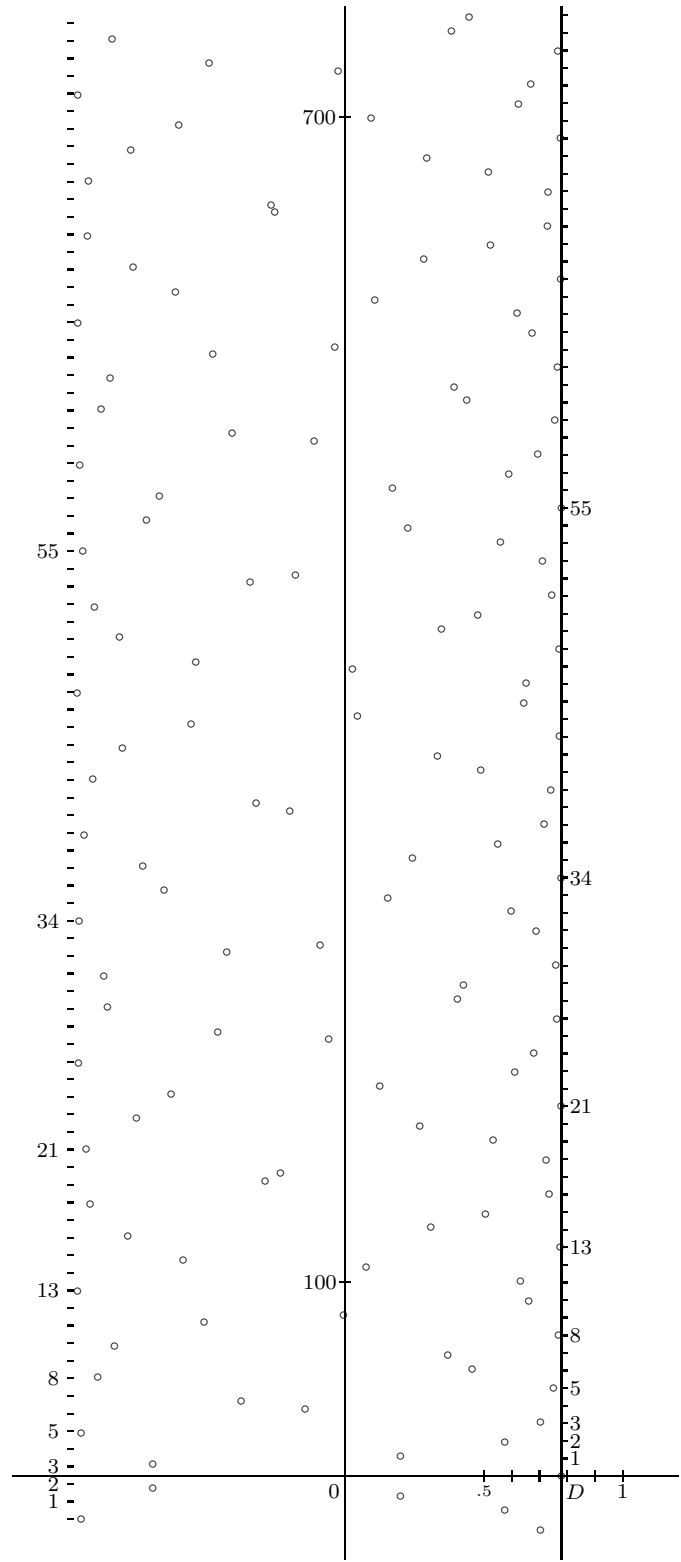
\begin{figure}[p]
\begin{picture}(324,580)(-155,-25)
\unitlength 1.05pt
\put(-120,0){\line(1,0){240}}\put(0,-30){\line(0,1){560}}\put(77.92122082,-30){\line(0,1){560}}
\put(-4,-5.5){\makebox(0,0){$\scriptstyle 0$}}
\put( 100   ,-2){\line(0,1){4}}\put( 100 ,-5.5){\makebox(0,0){$\scriptstyle 1$}}
\put( 90   ,-2){\line(0,1){4}}%\put( 90 ,-5.5){\makebox(0,0){$\scriptscriptstyle .9$}}
\put( 80   ,-2){\line(0,1){4}}%\put( 80 ,-5.5){\makebox(0,0){$\scriptscriptstyle .8$}}
\put( 70   ,-2){\line(0,1){4}}%\put( 70 ,-5.5){\makebox(0,0){$\scriptscriptstyle .7$}}
\put( 60   ,-2){\line(0,1){4}}%\put( 60 ,-5.5){\makebox(0,0){$\scriptscriptstyle .6$}}
\put( 50   ,-2){\line(0,1){4}}\put( 50 ,-5.5){\makebox(0,0){$\scriptscriptstyle .5$}}
\put( 83 ,-5.5){\makebox(0,0){$\scriptstyle D$}}
\put(-2, 70.0 ){\line(1,0){4}}\put(-9, 70.0 ){\makebox(0,0){$\scriptstyle 100$}}
\put(-2, 490.0 ){\line(1,0){4}}\put(-9.2, 490.0 ){\makebox(0,0){$\scriptstyle 700$}}

\put( 77.92122082 , 0 ){\makebox(0,0){$\nulpunt$}}	%D
%	aantal := 493	periode := 5529.479373	xschaal := 100	yschaal := .7
\put(-69.222 ,  -4.403){\makebox(0,0){$\nulpunt$}}
\put( 19.944 ,  -7.226){\makebox(0,0){$\nulpunt$}}
\put( 57.495 , -12.181){\makebox(0,0){$\nulpunt$}}
\put(-95.050 , -15.564){\makebox(0,0){$\nulpunt$}}
\put( 70.367 , -19.345){\makebox(0,0){$\nulpunt$}}

\put(-69.22260585 ,  4.4037070 ){\makebox(0,0){$\nulpunt$}}
\put( 19.94492484 ,  7.2269372 ){\makebox(0,0){$\nulpunt$}}
\put( 57.49503371 , 12.1810270 ){\makebox(0,0){$\nulpunt$}}
\put(-95.05035729 , 15.5646555 ){\makebox(0,0){$\nulpunt$}}
\put( 70.36771622 , 19.3456680 ){\makebox(0,0){$\nulpunt$}}
\put(-14.46493993 , 24.1878561 ){\makebox(0,0){$\nulpunt$}}
\put(-37.42913908 , 26.8249637 ){\makebox(0,0){$\nulpunt$}}
\put( 75.07168276 , 31.5361296 ){\makebox(0,0){$\nulpunt$}}
\put(-89.05890633 , 35.5632582 ){\makebox(0,0){$\nulpunt$}}
\put( 45.83446081 , 38.7129382 ){\makebox(0,0){$\nulpunt$}}
\put( 37.02549088 , 43.6894541 ){\makebox(0,0){$\nulpunt$}}
\put(-82.99170724 , 46.7043278 ){\makebox(0,0){$\nulpunt$}}
\put( 76.83731939 , 50.8799020 ){\makebox(0,0){$\nulpunt$}}
\put(-50.74463325 , 55.4857021 ){\makebox(0,0){$\nulpunt$}}
\put(  -.53695619 , 58.1701464 ){\makebox(0,0){$\nulpunt$}}
\put( 66.26037627 , 63.0682784 ){\makebox(0,0){$\nulpunt$}}
\put(-96.45095706 , 66.6963610 ){\makebox(0,0){$\nulpunt$}}
\put( 63.22901546 , 70.2302503 ){\makebox(0,0){$\nulpunt$}}
\put(  7.65764439 , 75.1543797 ){\makebox(0,0){$\nulpunt$}}
\put(-58.32270193 , 77.8839173 ){\makebox(0,0){$\nulpunt$}}
\put( 77.50846194 , 82.4164530 ){\makebox(0,0){$\nulpunt$}}
\put(-78.30573303 , 86.6814648 ){\makebox(0,0){$\nulpunt$}}
\put( 30.91341813 , 89.6184146 ){\makebox(0,0){$\nulpunt$}}
\put( 50.65787539 , 94.5897813 ){\makebox(0,0){$\nulpunt$}}
\put(-91.90113744 , 97.8270947 ){\makebox(0,0){$\nulpunt$}}
\put( 73.54803246 ,101.7606633 ){\makebox(0,0){$\nulpunt$}}
\put(-28.73555666 ,106.5279428 ){\makebox(0,0){$\nulpunt$}}
\put(-23.30285692 ,109.1534841 ){\makebox(0,0){$\nulpunt$}}
\put( 72.44276037 ,113.9516056 ){\makebox(0,0){$\nulpunt$}}
\put(-93.30915306 ,117.8271998 ){\makebox(0,0){$\nulpunt$}}
\put( 53.39800660 ,121.1196557 ){\makebox(0,0){$\nulpunt$}}
\put( 26.89482324 ,126.0858029 ){\makebox(0,0){$\nulpunt$}}
\put(-75.06535237 ,128.9774548 ){\makebox(0,0){$\nulpunt$}}
\put( 77.76314844 ,133.2962578 ){\makebox(0,0){$\nulpunt$}}
\put(-62.68064802 ,137.7810984 ){\makebox(0,0){$\nulpunt$}}
\put( 12.48655245 ,140.5435600 ){\makebox(0,0){$\nulpunt$}}
\put( 61.16950008 ,145.4807484 ){\makebox(0,0){$\nulpunt$}}
\put(-96.12188309 ,148.9570525 ){\makebox(0,0){$\nulpunt$}}
\put( 67.93768265 ,152.6433416 ){\makebox(0,0){$\nulpunt$}}
\put( -5.76167185 ,157.5222922 ){\makebox(0,0){$\nulpunt$}}
\put(-45.80876917 ,160.1847509 ){\makebox(0,0){$\nulpunt$}}
\put( 76.26474307 ,164.8327106 ){\makebox(0,0){$\nulpunt$}}
\put(-85.52034342 ,168.9522231 ){\makebox(0,0){$\nulpunt$}}
\put( 40.53357210 ,172.0173057 ){\makebox(0,0){$\nulpunt$}}
\put( 42.63019682 ,176.9950194 ){\makebox(0,0){$\nulpunt$}}
\put(-86.96448342 ,180.0922323 ){\makebox(0,0){$\nulpunt$}}
\put( 75.84457548 ,184.1763809 ){\makebox(0,0){$\nulpunt$}}
\put(-42.63299986 ,188.8494553 ){\makebox(0,0){$\nulpunt$}}
\put( -9.08321275 ,191.5005929 ){\makebox(0,0){$\nulpunt$}}
\put( 68.91509840 ,196.3662444 ){\makebox(0,0){$\nulpunt$}}
\put(-95.78859793 ,200.0889594 ){\makebox(0,0){$\nulpunt$}}
\put( 59.80500310 ,203.5294706 ){\makebox(0,0){$\nulpunt$}}
\put( 15.41376043 ,208.4737680 ){\makebox(0,0){$\nulpunt$}}
\put(-65.27766726 ,211.2586002 ){\makebox(0,0){$\nulpunt$}}
\put( 77.86137093 ,215.7127337 ){\makebox(0,0){$\nulpunt$}}
\put(-72.90229761 ,220.0645881 ){\makebox(0,0){$\nulpunt$}}
\put( 24.28308201 ,222.9291111 ){\makebox(0,0){$\nulpunt$}}
\put( 55.02507424 ,227.8911614 ){\makebox(0,0){$\nulpunt$}}
\put(-94.05859965 ,231.2184753 ){\makebox(0,0){$\nulpunt$}}
\put( 71.68594409 ,235.0577317 ){\makebox(0,0){$\nulpunt$}}
\put(-19.89817163 ,239.8736850 ){\makebox(0,0){$\nulpunt$}}
\put(-32.10590128 ,242.5019941 ){\makebox(0,0){$\nulpunt$}}
\put( 74.17363828 ,247.2485600 ){\makebox(0,0){$\nulpunt$}}
\put(-90.88962000 ,251.2182243 ){\makebox(0,0){$\nulpunt$}}
\put( 48.85787289 ,254.4215185 ){\makebox(0,0){$\nulpunt$}}
\put( 33.32514764 ,259.3953436 ){\makebox(0,0){$\nulpunt$}}
\put(-80.19287880 ,262.3616766 ){\makebox(0,0){$\nulpunt$}}
\put( 77.28813270 ,266.5925614 ){\makebox(0,0){$\nulpunt$}}
\put(-55.47770504 ,271.1539314 ){\makebox(0,0){$\nulpunt$}}
\put( 4.553845805 ,273.8648293 ){\makebox(0,0){$\nulpunt$}}
\put( 64.44047534 ,278.7796859 ){\makebox(0,0){$\nulpunt$}}
\put(-96.52926661 ,282.3497280 ){\makebox(0,0){$\nulpunt$}}
\put( 65.13986852 ,285.9415103 ){\makebox(0,0){$\nulpunt$}}
\put( 2.661993387 ,290.8503737 ){\makebox(0,0){$\nulpunt$}}
\put(-53.72811040 ,293.5508507 ){\makebox(0,0){$\nulpunt$}}
\put( 77.13328685 ,298.1291358 ){\makebox(0,0){$\nulpunt$}}
\put(-81.27421092 ,302.3393755 ){\makebox(0,0){$\nulpunt$}}
\put( 34.73552604 ,305.3237011 ){\makebox(0,0){$\nulpunt$}}
\put( 47.74497410 ,310.2987123 ){\makebox(0,0){$\nulpunt$}}
\put(-90.23382804 ,313.4817774 ){\makebox(0,0){$\nulpunt$}}
\put( 74.52517926 ,317.4730362 ){\makebox(0,0){$\nulpunt$}}
\put(-34.11367921 ,322.2066497 ){\makebox(0,0){$\nulpunt$}}
\put(-17.85708135 ,324.8376330 ){\makebox(0,0){$\nulpunt$}}
\put( 71.20719026 ,329.6637539 ){\makebox(0,0){$\nulpunt$}}
\put(-94.46034352 ,333.4811805 ){\makebox(0,0){$\nulpunt$}}
\put( 55.97035861 ,336.8295385 ){\makebox(0,0){$\nulpunt$}}
\put( 22.67915964 ,341.7888107 ){\makebox(0,0){$\nulpunt$}}
\put(-71.55383185 ,344.6375315 ){\makebox(0,0){$\nulpunt$}}
\put( 77.89784702 ,349.0089860 ){\makebox(0,0){$\nulpunt$}}
\put(-66.78768444 ,353.4443886 ){\makebox(0,0){$\nulpunt$}}
\put( 17.13522163 ,356.2432586 ){\makebox(0,0){$\nulpunt$}}
\put( 58.95721779 ,361.1914510 ){\makebox(0,0){$\nulpunt$}}
\put(-95.54137756 ,364.6105714 ){\makebox(0,0){$\nulpunt$}}
\put( 69.47575337 ,368.3551489 ){\makebox(0,0){$\nulpunt$}}
\put(-11.09299062 ,373.2123210 ){\makebox(0,0){$\nulpunt$}}
\put(-40.69790133 ,375.8575981 ){\makebox(0,0){$\nulpunt$}}
\put( 75.56930877 ,380.5452643 ){\makebox(0,0){$\nulpunt$}}
\put(-87.78096728 ,384.6080891 ){\makebox(0,0){$\nulpunt$}}
\put( 43.85329576 ,387.7248318 ){\makebox(0,0){$\nulpunt$}}
\put( 39.24034274 ,392.7022695 ){\makebox(0,0){$\nulpunt$}}
\put(-84.60363793 ,395.7483427 ){\makebox(0,0){$\nulpunt$}}
\put( 76.49378968 ,399.8889582 ){\makebox(0,0){$\nulpunt$}}
\put(-47.68184048 ,404.5213861 ){\makebox(0,0){$\nulpunt$}}
\put(-3.788327912 ,407.1915434 ){\makebox(0,0){$\nulpunt$}}
\put( 67.32518877 ,412.0779843 ){\makebox(0,0){$\nulpunt$}}
\put(-96.27470554 ,415.7424319 ){\makebox(0,0){$\nulpunt$}}
\put( 61.95949323 ,419.2402885 ){\makebox(0,0){$\nulpunt$}}
\put( 10.69621441 ,424.1728287 ){\makebox(0,0){$\nulpunt$}}
\put(-61.07515457 ,426.9225153 ){\makebox(0,0){$\nulpunt$}}
\put( 77.68172034 ,431.4254566 ){\makebox(0,0){$\nulpunt$}}
\put(-76.31234916 ,435.7242126 ){\makebox(0,0){$\nulpunt$}}
\put( 28.42489856 ,438.6325901 ){\makebox(0,0){$\nulpunt$}}
\put( 52.39038600 ,443.6008837 ){\makebox(0,0){$\nulpunt$}}
\put(-92.81136754 ,446.8725652 ){\makebox(0,0){$\nulpunt$}}
\put( 72.87238766 ,450.7699272 ){\makebox(0,0){$\nulpunt$}}
\put(-25.34388687 ,455.5568162 ){\makebox(0,0){$\nulpunt$}}
\put(-26.70246487 ,458.1817175 ){\makebox(0,0){$\nulpunt$}}
\put( 73.14867691 ,462.9608959 ){\makebox(0,0){$\nulpunt$}}
\put(-92.45938973 ,466.8727524 ){\makebox(0,0){$\nulpunt$}}
\put( 51.70540205 ,470.1306468 ){\makebox(0,0){$\nulpunt$}}
\put( 29.42953757 ,475.1002446 ){\makebox(0,0){$\nulpunt$}}
\put(-77.12245350 ,478.0199493 ){\makebox(0,0){$\nulpunt$}}
\put( 77.61805610 ,482.3052560 ){\makebox(0,0){$\nulpunt$}}
\put(-59.98552951 ,486.8201231 ){\makebox(0,0){$\nulpunt$}}
\put(  9.48885626 ,489.5615749 ){\makebox(0,0){$\nulpunt$}}
\put( 62.47433373 ,494.4908518 ){\makebox(0,0){$\nulpunt$}}
\put(-96.35697038 ,498.0031056 ){\makebox(0,0){$\nulpunt$}}
\put( 66.90588724 ,501.6530008 ){\makebox(0,0){$\nulpunt$}}
\put( -2.48188560 ,506.5442357 ){\makebox(0,0){$\nulpunt$}}
\put(-48.91607020 ,509.2198927 ){\makebox(0,0){$\nulpunt$}}
\put( 76.63692723 ,513.8417789 ){\makebox(0,0){$\nulpunt$}}
\put(-83.97149676 ,517.9964111 ){\makebox(0,0){$\nulpunt$}}
\put( 38.36335944 ,521.0299298 ){\makebox(0,0){$\nulpunt$}}
\put( 44.65431380 ,526.0070682 ){\makebox(0,0){$\nulpunt$}}
%	~150 puntjes

\put( 77.92122082 , 0 ){
\put(3,   6.3454   ){{\makebox(0,0)[l]{$\scriptstyle   1   $}}}
\put(3,   12.691   ){{\makebox(0,0)[l]{$\scriptstyle   2   $}}}
\put(3,   19.036   ){{\makebox(0,0)[l]{$\scriptstyle   3   $}}}
\put(3,   31.727   ){{\makebox(0,0)[l]{$\scriptstyle   5   $}}}
\put(3,   50.763   ){{\makebox(0,0)[l]{$\scriptstyle   8   $}}}
\put(3,   82.490   ){{\makebox(0,0)[l]{$\scriptstyle   13   $}}}
\put(3,   133.25   ){{\makebox(0,0)[l]{$\scriptstyle   21   $}}}
\put(3,   215.74   ){{\makebox(0,0)[l]{$\scriptstyle   34   $}}}
\put(3,   348.99   ){{\makebox(0,0)[l]{$\scriptstyle   55   $}}}
\multiput(0,0)(0,6.3454){84}{\line(1,0){2}}}
\put(-100.050 , -15.564){
\put(-3,   6.3454   ){{\makebox(0,0)[r]{$\scriptstyle   1   $}}}
\put(-3,   12.691   ){{\makebox(0,0)[r]{$\scriptstyle   2   $}}}
\put(-3,   19.036   ){{\makebox(0,0)[r]{$\scriptstyle   3   $}}}
\put(-3,   31.727   ){{\makebox(0,0)[r]{$\scriptstyle   5   $}}}
\put(-3,   50.763   ){{\makebox(0,0)[r]{$\scriptstyle   8   $}}}
\put(-3,   82.490   ){{\makebox(0,0)[r]{$\scriptstyle   13   $}}}
\put(-3,   133.25   ){{\makebox(0,0)[r]{$\scriptstyle   21   $}}}
\put(-3,   215.74   ){{\makebox(0,0)[r]{$\scriptstyle   34   $}}}
\put(-3,   348.99   ){{\makebox(0,0)[r]{$\scriptstyle   55   $}}}
\multiput(0,0)(0,6.3454){86}{\line(1,0){2}}}
\end{picture}
\caption[The almost periodic behavior of the complex dimensions of the golden flow]
{The almost periodic behavior of the complex dimensions of the golden flow.}
\label{F: golden global}
\end{figure}

A diagram of the complex dimensions of the golden flow is given in Figure~\ref{F: golden global}.
To obtain it,
we chose the approximation $\phi\approx 987/610$ to approximate the flow by the lattice flow with weights
$w_1=610 w$,
$w_2=987 w$,
where $w=(1/610)\log 2$.
We then used Maple to solve the corresponding polynomial equation.
In particular,
the dimension~$D$ of the golden flow is approximately equal to $D=.77921\dots$.
See also Example~\ref{E: cxd of GF}.
\end{example}

\section{Explicit Formulas}
\label{S: explicit}

We will formulate our explicit formulas in the more general framework of~\cite[Chapter~4]{long}.
We refer to this book for the proofs and much additional information.

Let $\GS$ be a positive measure on $(0,\infty )$,
supported away from $0$.
Its Mellin transform is
\begin{gather}\label{E: zeta eta}
\zeta_\GS (s)=\int_0^\infty x^{-s}d\GS ,
\end{gather}
the {\em geometric zeta function\/} of $\GS$.
We assume that $\zeta_\GS$ is convergent for some $s$,
and we write $D$ for the abcissa of convergence.
We assume that there exists a screen to the left of $\Re s=D$ such that~$\zeta_\GS$ has a meromorphic
continuation to the corresponding window.
To simplify the exposition,
we also assume that the poles of $\zeta_\eta$ are simple.
This is the case,
for example,
if $\zeta_\weight$ has a meromorphic continuation to a neighborhood of $W$.
The general case,
when the poles of $\zeta_\GS$ may have arbitrary multiplicities,
is treated in~\cite[Chapter~4]{long}.

The screen $S$ is given as the graph of a bounded function $S$,
with the horizontal and vertical axes interchanged:
$$
S=\{ S(t)+it\colon t\in\Real\} .
$$
We will write $\inf S=\inf_{t\in\Real}S(t)$ and $\sup S=\sup_{t\in\Real}S(t).$
We assume in addition that $S$ is a Lipschitz continuous function;
i.e.,
there exists a nonnegative real number,
denoted by $\| S\|_{\rm Lip}$,
such that
$$
|S(x)-S(y)|\leq\| S\|_{\rm Lip}|x-y|\qquad\mbox{for all }x,y\in\Real .
$$
Further,
recall from Section~\ref{S: DS} that the window $W$ is the part of the complex plane to the
right of $S$;
see formula~(\ref{E: window}).
\medskip

Assume that $\zeta_\GS$ satisfies the following growth conditions:
\begin{quote}
There exist real constants $\kappa\geq 0$ and $C>0$ and a sequence $\{ T_n\}_{n\in\Int}$ of real numbers
tending to $\pm\infty$ as $n\to\pm\infty$,
with $T_{-n}<0<T_n$ for $n\geq 1$ and $\lim_{n\to +\infty}T_n/|T_{-n}|=1$,
such that\medskip
\begin{description}
\item[(H$_1$)]
For all $n\in\Int$ and all $\sigma\geq S(T_n)$,
\begin{align}\label{E: H1}
|\zeta_\GS (\sigma +iT_n)|\leq C\cdot |T_n|^\kappa ,
\end{align}
\item[(H$_2$)]
For all $t\in\Real$,
$|t|\geq 1$,
\begin{align}\label{E: H2}
|\zeta_\GS (S(t)+it)|\leq C\cdot |t|^\kappa .
\end{align}
\end{description}
\end{quote}

Hypothesis ({\bf H$_1$}) is a polynomial growth condition along horizontal lines
(necessarily avoiding the poles of $\zeta_\GS$),
while hypothesis ({\bf H$_2$}) is a polynomial growth condition along the vertical direction of the
screen.

In the following,
we denote by $\res (g(s);\omega )$ the residue of a meromorphic function $g=g(s)$ at $\omega$.
It vanishes unless $\omega$ is a pole of $g$.
Also,
\begin{align}
(s)_k=\frac{\Gamma (s+k)}{\Gamma (s)},
\end{align}
for $k\in\Int$.
Thus,
$(s)_0=1$ and,
for $k\geq 1$,
$(s)_k=s(s+1)\ldots (s+k-1)$.
\medskip

Let
\begin{gather}\label{E: Neta}
N_\GS (x)=N^{[1]}_\GS (x)=\GS (0,x)+\frac{1}{2}\GS (\{ x\} ),
\end{gather}
and more generally,
let $N^{[k]}_\GS (x)$ be the $(k-1)$-st antiderivative of this function,
for $k=1,2,\dots$.
Our explicit formula expresses this function as a sum over the poles of $\zeta_\GS$.

\begin{theorem}[The pointwise explicit formula]
\label{T: pointwise}
Let $\GS$ be a generalized fractal string,
satisfying hypotheses\/ {\rm ({\bf H$_1$})} and\/ {\rm ({\bf H$_2$}).}
Let\/ $k$ be a positive integer such that\/ $k>\kappa +1,$
where $\kappa\geq 0$ is the exponent occurring in the statement of\/ {\rm (\bf H$_1$\rm )} and\/
{\rm ({\bf H$_2$})}.
Then,
for all\/ $x>0,$
the pointwise explicit formula is given by the following equality\/{\rm :}
\begin{equation}\label{E: pointwise with error}
\begin{aligned}
N^{[k]}_\GS (x)&
=\sum_{\substack{\omega\in\Dimensions_\GS (W)\\ \omega\not\in\{ 0,-1,\dots ,-(k-1)\}}}
\res\left(\zeta_\GS (s);\omega\right)\frac{x^{\omega +k-1}}{(\omega )_k}\\
&\hspace{5mm}
+\sum_{j=0}^{k-1}\res\left(\frac{x^{s+k-1}\zeta_\GS (s)}{(s)_k};-j\right)
+R^{[k]}_\GS (x).
\end{aligned}
\end{equation}
Here,
for $x>0,$
$R(x)=R^{[k]}_\GS (x)$ is the error term,
given by the absolutely convergent integral
\begin{gather}\label{E: pointwise error formula}
R(x)=R^{[k]}_\GS (x)=\frac{1}{2\pi i}\int_{S}x^{s+k-1}\zeta_\GS (s)\frac{ds}{(s)_k} .
\end{gather}
Further,
for all $x>0,$
we have
\begin{align}\label{E: estimate of error}
R(x)=R^{[k]}_\GS (x)\leq C(1+\| r\|_{\rm Lip})\frac{x^{k-1}}{k-\kappa -1}
\max\{ x^{\sup S}, x^{\inf S}\} +C',
\end{align}
where $C$ is the positive constant occurring in\/ {\rm ({\bf H$_1$})} and\/ {\rm ({\bf H$_2$})} and $C'$ is
some suitable positive constant.
The constants $C(1+\| r\|_{\rm Lip})$ and $C'$ depend only on $\GS$ and the screen,
but not on $k$.

In particular,
we have the following pointwise error estimate{\rm :}
\begin{gather}\label{E: order of error}
R(x)=R^{[k]}_\GS (x)=O\left( x^{\sup S+k-1}\right) ,
\end{gather}
as $x\to\infty$.
Moreover,
if\/ $S(t)<\sup S$ for all\/ $t\in\Real$ {\rm (}i.e.,
if the screen lies strictly to the left of the line $\Re s=\sup S),$
then $R(x)$ is of order less than $x^{\sup S+k-1}$ as $x\to\infty\colon$
\begin{gather}\label{E: pointwise R=o}
R(x)=R^{[k]}_\GS (x)=o\left( x^{\sup S+k-1}\right) ,
\end{gather}
as $x\to\infty$.
\end{theorem}

To formulate our second explicit formula,
a distributional formula,
we view $\GS$ as a distribution,
acting on a test function $\varphi$ defined on $(0,\infty )$ by
\begin{gather}\label{E: eta}
\left< N^{[0]}_\GS ,\varphi\right> =\int_0^\infty\varphi\,d\GS .
\end{gather}
We then define $N^{[k]}_\GS$ as the distribution obtained by integrating this one $k$ times,
so that
\begin{gather}\label{E: eta k}
\left< N^{[k]}_\GS ,\varphi\right>
=\int_0^\infty\int_y^\infty\frac{(x-y)^{k-1}}{(k-1)!}\varphi (x)\, dx\,\GS (dy).
\end{gather}
It is easily verified that this definition coincides with formula~(\ref{E: Neta}) and the next line above when
$k\geq 1$.
For $k\leq 0$,
we extend this definition by differentiating $|k|$ times the distribution $N^{[0]}_\GS$.

We shall denote by $\widetilde{\varphi}$ the {\em Mellin transform\/}
of a (suitable) function~$\varphi$ on $(0,\infty )$;
it is defined by
\begin{equation}\label{E: Mellin}
\widetilde{\varphi}(s)=\int_0^\infty\varphi (x)x^{s-1}\, dx.
\end{equation}

\begin{theorem}[The distributional explicit formula]
\label{T: distributional}
Let $\GS$ be a generalized fractal string satisfying hypotheses {\rm ({\bf H$_1$})} and\/
{\rm (\bf H$_2$\rm ).}
Then,
for every $k\in\Int ,$
the distribution $\primitive^{[k]}_\GS$ is given by formula\/~{\rm (\ref{E: pointwise with error}),}
interpreted as a distribution.
That is,
the action of $\primitive^{[k]}_\GS$ on a test function $\varphi$ is given by
\begin{equation}\label{E: distributional with error test}
\begin{aligned}
\left<\primitive^{[k]}_\GS ,\varphi\right> &=
\sum_{\substack{\omega\in\Dimensions_\GS (W)\\ \omega\not\in\{ 0,-1,\dots ,-(k-1)\}}}
\res\left(\zeta_\GS (s);\omega\right)\frac{\widetilde{\varphi}(\omega +k)}{(\omega )_k}
\\&\hspace{5mm}
+\sum_{j=0}^{k-1}\res\left(\frac{\zeta_\GS (s)\widetilde{\varphi}(s+k)}{(s)_k};-j\right)
+\left< \Rerror^{[k]}_\GS ,\varphi\right> .
\end{aligned}
\end{equation}
Here,
the distribution $\Rerror =\Rerror^{[k]}_\GS$ is the error term,
given by
\begin{gather}\label{E: distributional error}
\left< \Rerror,\varphi\right> =\left< \Rerror_\GS^{[k]},\varphi\right> =
\frac{1}{2\pi i}\int_{S}\zeta_\GS (s)\widetilde{\varphi}(s+k)\frac{ds}{(s)_k}.
\end{gather}
\end{theorem}

\begin{definition}\label{D: shifted dirac, order of distribution}
We will say that a distribution $\Rerror$ on $(0,\infty )$ is of {\em asymptotic order\/} at most $x^\alpha$
(respectively,
less than $x^\alpha$)---and we will write $\Rerror (x)=O(x^\alpha )$ (respectively,
$\Rerror (x)=o(x^\alpha )$),
as $x\to\infty$---if applied to a test function $\varphi$,
we have that
\begin{align}
\left<\Rerror ,\varphi_a\right> =O\left( a^\alpha\right) \quad\mbox{(respectively, }
\left<\Rerror ,\varphi_a\right> =o\left( a^\alpha\right) ),
\quad\mbox{as }a\to\infty ,
\end{align}
where $\varphi_a(x)=a^{-1}\varphi (x/a)$.
\end{definition}

\begin{theorem}\label{T: error estimate}
Fix $k\in\Int$.
Assume that the hypotheses of Theorem~\ref{T: distributional} are satisfied,
and let the distribution $\Rerror =\Rerror^{[k]}_\GS$ be given by\/~{\rm (\ref{E: distributional error}).}
Then $\Rerror$ is of asymptotic order at most $x^{\sup S+k-1}$ as $x\to\infty\colon$
\begin{align}\label{E: R=O}
\Rerror^{[k]}_\GS (x)=O\left( x^{\sup S+k-1}\right) ,
\quad\mbox{as }x\to\infty ,
\end{align}
in the sense of Definition~\ref{D: shifted dirac, order of distribution}.

Moreover,
if\/ $S(t)<\sup S$ for all\/ $t\in\Real$ {\rm (}i.e.,
if the screen lies strictly to the left of the line $\Re s=\sup S),$
then $\Rerror$ is of asymptotic order less than $x^{\sup S+k-1}$ as $x\to\infty\colon$
\begin{align}\label{E: R=o}
\Rerror^{[k]}_\GS (x)=o\left( x^{\sup S+k-1}\right) ,
\quad\mbox{as }x\to\infty .
\end{align}
\end{theorem}

We refer to \cite[Chapter~4]{long} for a proof of Theorems~\ref{T: pointwise},
\ref{T: distributional} and~\ref{T: error estimate}.

\begin{remark}[Oscillatory terms in the explicit formula]
Our explicit formulas give expansions of various functions associated with a measure $\GS$ as a sum over the
poles of $\zeta_\GS$.
The term corresponding to the pole $\omega$ of multiplicity one is of the form $Cx^\omega$,
where $C$ is a constant depending on $\omega$.
If $\omega$ is real,
the function $x^\omega$ simply has a certain asymptotic behavior as $x\to\infty$.
If,
on the other hand,
$\omega =\beta +i\gamma$ has a nonzero imaginary part $\gamma$,
then $x^\omega =x^\beta\cdot x^{i\gamma}$ is of order $O(x^\beta )$ as $x\to\infty$,
with a multiplicatively periodic behavior:
The function $x^{i\gamma}=\exp (i\gamma\log x)$ takes the same value at the points $e^{2\pi n/\gamma}x$
($n\in\Int$).
Thus,
the term corresponding to $\omega$ will be called an oscillatory term.
If there are poles with higher multiplicity,
there will also be terms of the form $Cx^\omega (\log x)^m$,
$m\in\Nat^*$,
which have a similar oscillatory behavior.
\end{remark}

\begin{example}[The classical Prime Number Theorem]\label{E: Riemann-von Mangoldt}
Let $\zeta (s)=1+2^{-s}+3^{-s}+\dots$ (for $\Re s>1$) be the Riemann zeta function.
This function has an Euler product $\zeta (s)=\prod_p(1-p^{-s})^{-1}$ (for $\Re s>1$),
where $p$ runs over the prime numbers.
Analogously to Corollary~\ref{C: logarithmic derivative and psi},
we obtain $-\zeta'/\zeta (s)=\int_0^\infty x^{-s}\, d\psi (x)=\zeta_\primestring (s)$,
where $\psi (x)=\sum_{p^k\leq x}\log p=N^{[1]}_\primestring (x)$,
and
\begin{align}\label{E: prime string}
\primestring =\sum_{m\geq 1,\, p}(\log p)\delta_{\{ p^m\}}
\end{align}
is the {\em prime string\/} (see~\cite{long}).
We apply Theorem~\ref{T: distributional} to $\GS =\primestring$ to obtain the explicit formula for $\psi$:
\begin{gather}\label{E: psi}
\psi (x)=x-\sum_{\rho\in W}\res (\zeta'/\zeta (s);\rho )\frac{x^\rho}{\rho}
-\frac{1}{2\pi i}\int_S\frac{\zeta'}{\zeta}(s)x^s\frac{ds}{s},
\end{gather}
where $\rho$ runs through the sequence of critical zeros of $\zeta$:
$\zeta (\rho )=0$,
$0<\Re\rho <1$.

By means of classical arguments~\cite[Theorem~19]{Ingham},
it is known that $\zeta$ has a zero free region of the form
$$
\left\{\sigma +it\in\Com\colon \sigma >1-C/\log t\right\} ,
$$
for some positive constant $C$.
Also,
$-\zeta'/\zeta$ is not too large on the boundary of this region.
In our language,
this means that we can choose a screen to the left of $\Re s=1$ such that there are no zeros of $\zeta$ in $W$
and $-\zeta'/\zeta$ satisfies ({\bf H$_1$}) and ({\bf H$_2$}).
Then~(\ref{E: psi}) becomes
\begin{gather}
\psi (x)=x+o(x),
\end{gather}
as $x\to\infty$.
This is equivalent to the classical Prime Number Theorem.

Using the existence of the zero free region,
one can derive the stronger estimate
\begin{gather}\label{E: classical PNT}
\psi (x)=x+O\left( xe^{-c\sqrt{\log x}}\right) ,
\end{gather}
as $x\to\infty$,
for some positive constant $c$ (see~[\citen{hE74}; \citen{Ingham}, Theorem~23]).
This is the classical Prime Number Theorem,
with Error Term.
\end{example}

\section{The Prime Orbit Theorem for Flows}
\label{S: PNT}

Let $\Flow_\weight$ be a suspended flow as in Section~\ref{S: DS}.
In Corollary~\ref{C: logarithmic derivative and psi},
we have written the logarithmic derivative of $\zeta_\weight (s)$ as the Mellin transform of the counting
function $\psi_\weight$ of the weighted periodic orbits of $\sigma$,
as defined in~(\ref{E: psi w}).
Put $\GS =d\psi_\weight$,
so that $\zeta_\GS =-\zeta_\weight'/\zeta_\weight$.
The poles of $-\zeta_\weight'/\zeta_\weight$ are the complex dimensions of $\Flow_\weight$ and
the residue at $\omega$ is $-\ord (\zeta_\weight ;\omega )$.
By Theorems~\ref{T: distributional} and~\ref{T: error estimate},
we obtain the following explicit formula for the counting function of weighted periodic orbits of $\sigma$.

\begin{theorem}[The Prime Orbit Theorem with Error Term]\label{T: PNT}
Let\/ $\Flow_\weight$ be a suspended flow that satisfies conditions {\rm ({\bf H$_1$})} and
{\rm (\bf H$_2$\rm ).}
Then we have the following equality between distributions\/{\rm :}
\begin{gather}\label{E: PNT}
\psi_\weight (x)=\frac{x^D}{D}
+\sum_{\omega\in\Dimensions_\weight\backslash\{ D,0\}}
-\ord\left(\zeta_\weight ;\omega\right)\frac{x^\omega}{\omega}
+\res\left( -\frac{x^s\zeta_\weight'(s)}{s\zeta_\weight (s)};0\right) +R (x),
\end{gather}
where $\ord\left(\zeta_\weight ;\omega\right) <0$ denotes the order of\/ $\zeta_\weight$ at\/ $\omega ,$
and
\begin{gather}
R (x)=-\int_S\frac{\zeta_\weight'}{\zeta_\weight}(s)x^s\, \frac{ds}{s}
=O\left( x^{\sup S}\right) ,
\end{gather}
as $x\to\infty .$

If\/ $0$ is not a complex dimension of the flow,
then the third term on the right-hand side of\/ {\rm (\ref{E: PNT})} simplifies to
$
-\zeta_\weight'/\zeta_\weight (0).
$
In general,
this term is of the form $p+q\log x,$
for some constants $p$ and\/ $q.$

If\/ $D$ is the only complex dimension on the line $\Re s=D,$
then the error term,
\begin{gather}\label{E: et}
\sum_{\omega\in\Dimensions_\weight\backslash\{ D,0\}}
-\ord\left(\zeta_\weight ;\omega\right)\frac{x^\omega}{\omega}
+\res\left( -\frac{x^s\zeta_\weight'(s)}{s\zeta_\weight (s)};0\right) +R (x),
\end{gather}
is estimated by $o(x^D),$
as $x\to\infty .$
If this is the case,
then we obtain a Prime Orbit Theorem for $\Flow_\weight$ as follows:
\begin{gather}
\psi_\weight (x)=\frac{x^D}{D}+o\left( x^D\right) ,
\end{gather}
as $x\to\infty .$
\end{theorem}

\begin{proof}
The first part of the theorem follows from the distributional explicit formula
(Theorem~\ref{T: distributional}) and from the first part of Theorem~\ref{T: error estimate},
while the second part follows from the second part of Theorem~\ref{T: error estimate}.
\end{proof}

The explicit formula holds for every flow satisfying our conditions on $-\zeta_\weight'/\zeta_\weight$.
In particular,
we can apply it to the `axiom A flows' considered in~\cite[Chapter~6]{wPmP90},
in view of~\cite[pp.~100--101]{wPmP90}.
We hope to do so more explicitly in a later work.
For simplicity and for the sake of concision,
however,
we will focus in the rest of this paper on the important example of `self-similar flows' (in the sense of
Section~\ref{S: ssd} above).

\section{The Prime Orbit Theorem for Self-Similar Flows}
\label{S: PNT SS}

For self-similar flows,
$\zeta_\weight$ does not have any zeros (see~(\ref{E: dynamical zeta})).
Hence every contribution to (\ref{E: PNT}) comes from a pole of $\zeta_\weight$,
and each coefficient $-\ord (\zeta_\weight ;\omega )$ is positive.
Furthermore,
$0$ is never a complex dimension,
so the third term on the right-hand side of\/ {\rm (\ref{E: PNT})} in the explicit formula is
\begin{gather}
-\frac{\zeta_\weight'(0)}{\zeta_\weight (0)}=-\frac{1}{N-1}\sum_{j=1}^Nw_j.
\end{gather}
We can obtain information about $\psi_\weight$ by choosing a suitable screen.

\subsection{Lattice Flows}
\label{S: geo lattice}

In the lattice case,
we obtain the Prime Orbit Theorem for lattice self-similar flows:
\begin{gather}\label{E: PNT psi lattice}
\psi_\weight (x)=g_1(\log x)x^D-\frac{1}{N-1}\sum_{j=1}^Nw_j+O\left( x^{D-\alpha}\right) ,
\end{gather}
as $x\to\infty$.
Here,
$D-\alpha$ is the abcissa of the first vertical line of complex dimensions next to $D$,
and the periodic function~$g_1$,
of period $w$,
is given by\footnote
{We use the notation $\{ u\}$ for the fractional part,
and $[u]$ for the integer part of the real number $u$,
so that $\{ u\} =u-[u]\in [0,1)$.}
\begin{gather}\label{E: PNT psi g1}
g_1(y)=\sum_{n=-\infty}^\infty\frac{e^{2\pi iny/w}}{D+2\pi in/w}
=\frac{b_1w}{b_1-1}b_1^{-\{ y/w\}},
\end{gather}
where $b_1=e^{wD}$.
By choosing a screen located to the left of all the complex dimensions of $\Flow_\weight$,
we can even obtain more precise information about $\psi_\weight$.
In the notation of Theorem~\ref{T: struc cxd},
we obtain
\begin{equation}\label{E: PNT psi lattice precise}
\begin{aligned}
\psi_\weight (x)&=\sum_{u=1}^q -\ord\left(\zeta_\weight ;\omega_u\right)
\sum_{n\in\Int}\frac{x^{\omega_u+2\pi in/w}}{\omega_u+2\pi in/w}-\frac{1}{N-1}\sum_{j=1}^Nw_j\\
&=\sum_{u=1}^q -\ord\left(\zeta_\weight ;\omega_u\right) g_u(\log x)x^{\omega_u}-\frac{1}{N-1}\sum_{j=1}^Nw_j,
\end{aligned}
\end{equation}
where for each $u=1,\dots ,q$,
the function $g_u$ is periodic of period $w$,
given by
\begin{gather}
g_u(y)=\sum_{n\in\Int}\frac{e^{2\pi iny/w}}{\omega_u+2\pi in/w}=\frac{b_uw}{b_u-1}b_u^{-\{ y/w\}},
\end{gather}
where $b_u=e^{w\omega_u}$.
Here,
$\omega_1(=D),\omega_2,\dots ,\omega_q$ are given as in the lattice case of Theorem~\ref{T: struc cxd} and
$\ord (\zeta_\weight ;\omega_1 )=-1$.

For instance,
the Cantor flow (with $D=\log_3 2$ and $w_1=w_2=w=\log 3$,
see Example~\ref{E: Cantor Flow}) has
\begin{gather}
\psi_{\rm CF} (x)=g_1(\log x)x^D-2\log 3,
\end{gather}
with $g_1(y)=w2^{1-\{ y/w\}}$,
while the Fibonacci flow\footnote
{Also called the golden mean flow in the literature (see,
e.g.,
\cite[p.~59]{tBmKcS91}),
but not to be confused with the golden flow in our present paper,
which is a nonlattice self-similar flow.}
of Example~\ref{E: Fibonacci Flow} (with $D=\log_2\phi$ and $w_1=w=\log 2$,
$w_2=2w$) has:
\begin{gather}
\psi_{\rm Fib} (x)=g_1(\log x)x^D+g_2(\log x)x^{\pi i/w}x^{-D}-3\log 2,
\end{gather}
where $g_1(y)=w\phi^{2-\{ y/w\}}$ and
$$
g_2(y)=\sum_{n\in\Int}\frac{e^{2\pi iny/w}}{-D+2\pi i(n+1/2)/w}=w\phi^{\{y/w\} -2}e^{-\pi i\{ y/w\}}.
$$
In the second asymptotic term,
the product $e^{-\pi i\{ (\log x)/w\}}x^{\pi i/w}$ combines to give the sign $(-1)^{[(\log x)/w]}$.

\subsection{Nonlattice Flows}
\label{S: geo nonlattice}

In the nonlattice case,
we use Theorem~\ref{T: struc cxd} according to which there exists $\delta >0$ and a screen $S$ lying to the
left of the vertical line $\Re s=D-\delta$ such that $-\zeta_\weight'/\zeta_\weight$ is bounded on $S$ and all
the complex dimensions $\omega$ to the right of $S$ have residue
$\res (-\zeta_\weight'/\zeta_\weight ;\omega )=1$.
Then $R (x)=O(x^{D-\delta})$,
as $x\to\infty$.
There are no complex dimensions with $\Re\omega =D$ except for $D$ itself.
Hence,
the assumptions of Theorem~\ref{T: PNT} are satisfied.
Therefore,
in view of Theorem~\ref{T: PNT},
we deduce by a classical argument (see the proof of Theorems~\ref{T: nl PNT} and~\ref{T: nl PNT N} on
page~\pageref{P: nl PNT N}) the Prime Orbit Theorem for nonlattice suspended flows:
\begin{gather}\label{E: PNT psi nonlattice}
\psi_\weight (x)=\frac{x^D}{D}
+\sum_{\omega\in\Dimensions_\weight\backslash\{ D\}}\frac{x^\omega}{\omega}+O(x^{D-\delta})
=\frac{x^D}{D}+o\left( x^D\right) ,
\end{gather}
as $x\to\infty$.
(See also Theorem~\ref{T: nl PNT N},
and when $N=2$,
Theorem~\ref{T: nl PNT} below for a better estimate of the error.)
We note that this estimate is always best possible,
since by the nonlattice case of Theorem~\ref{T: struc cxd},
there always exist complex dimensions of $\weight$ arbitrarily close to the vertical line $\Re s=D$.

\begin{remark}
It would be interesting to apply Theorem~\ref{T: PNT} to suspended flows that are more general than
self-similar flows:
for example,
those considered by Lalley in~\cite[2]{sL89},
such as the `approximately self-similar flows' naturally associated with limit sets of suitable Kleinian
groups.
This would require a more detailed study of the dynamical zeta function of each of these flows.
It is known that the lattice-nonlattice dichotomy applies in these more general cases;
see~\cite[2]{sL89}.
We hope to investigate this question in a later work.
\end{remark}

\section{The Error Term in the Nonlattice Case}
\label{S: error nl}

A nonlattice flow has weights $w_1\leq\dots\leq w_N$,
where at least one ratio $w_j/w_k$ is irrational.
Let
\begin{gather}\label{E: f}
f(s)=1-\sum_{j=1}^Ne^{-w_js}.
\end{gather}
Then $D>0$ is the unique real solution of the equation $f(s)=0$.
Moreover,
the derivative
\begin{gather}
f'(s)=\sum_{j=1}^Nw_je^{-w_js}
\end{gather}
does not vanish at $D$.

When $N=2$,
the flow is called a Bernoulli flow (see Section~\ref{S: ssd}).
We write $\alpha =w_2/w_1$,
so $\alpha >1$ is irrational.
In this case,
we obtain very detailed information about the growth of $-\zeta_\weight'/\zeta_\weight$ on the line $\Re s=D$,
and we can compute a pole free region for this function,
by considering the continued fraction expansion of $\alpha$.
We briefly collect here the facts that we will use.
See~\cite{HW,gelijkverdeling},
and~\cite{lB,thesis,subspaces} for a connection with the Riemann Hypothesis.

\subsection{Continued Fractions}

Let $\alpha$ be an irrational real number with a continued fraction expansion
$\alpha =[[a_0,a_1,a_2,\dots ]]=a_0+1/(a_1+1/(a_2+\dots ))$.
We recall that the two sequences $a_0,a_1,\dots$ and $\alpha_0,\alpha_1,\dots$ are defined by
$\alpha_0=\alpha$ and,
for $n\geq 0$,
$a_n=[\alpha_n]$,
the integer part of $\alpha_n$,
and $\alpha_{n+1}=1/({\alpha_n-a_n})$.
The {\em convergents\/} of $\alpha$,
\begin{gather}\label{E: p/q}
\frac{p_n}{q_n}=[[a_0,a_1,a_2,\dots ,a_n]],
\end{gather}
are successively computed by $p_{-2}=q_{-1}=0$,
$p_{-1}=q_{-2}=1$,
and $p_{n+1}=a_{n+1}p_n+p_{n-1}$,
$q_{n+1}=a_{n+1}q_n+q_{n-1}$.
We also define $q'_n=\alpha_1\cdot\alpha_2\cdot\dots\cdot\alpha_n$,
and note the formula $q'_{n+1}=\alpha_{n+1}q_n+q_{n-1}$.
Then
\begin{align}\label{approx}
q_n\alpha -p_n=\frac{(-1)^n}{q'_{n+1}}.
\end{align}
We have $q_n\geq\phi^{n-1}$,
where $\phi =({1+\sqrt{5}})/2$ is the golden ratio.

Let $n\in\Nat$ and choose $l$ such that $q_{l+1}> n$.
We can successively compute (see~\cite{gelijkverdeling})
$$
n=d_lq_l+n_l,\ n_l=d_{l-1}q_{l-1}+n_{l-1},\dots ,n_1=d_0q_0,
$$
where $d_\nu$ is the quotient and $n_\nu<q_\nu$ is the remainder of the division of $n_{\nu +1}$ by~$q_\nu$.
We set $d_{l+1}=d_{l+2}=\ldots =0$.
Then
\begin{align}\label{E: adic expansion}
n=\sum_{\nu =0}^\infty d_\nu q_\nu .
\end{align}
We call this the {\em$\alpha$-adic expansion\/} of $n$.
Note that $0\leq d_\nu\leq a_{\nu +1}$ and that if $d_\nu =a_{\nu +1}$,
then $d_{\nu -1}=0$.
Also $d_0<a_1$.
It is not difficult to show that these properties uniquely determine the sequence $d_0,d_1,\dots$ of
{\em $\alpha$-adic digits\/} of $\alpha$.

\begin{lemma}\label{L: adic expansion}
Let\/ $n$ be given by\/~{\rm (\ref{E: adic expansion}).}
Let\/ $k\geq 0$ be such that\/ $d_k\neq 0$ and\/ $d_{k-1}=\dots =d_0=0.$
Put\/ $m=\sum_{\nu =k}^\infty d_\nu p_\nu .$
Then $n\alpha -m$ lies strictly between $(-1)^k/q'_{k+2}$ and $(-1)^k/q'_k.$
\end{lemma}

\begin{proof}
We have $n\alpha -m=\sum_{\nu =k}^\infty d_\nu (\alpha q_\nu -p_\nu ).$
Since $\alpha q_\nu -p_\nu =(-1)^\nu /q'_{\nu +1}$,
the terms in this sum are alternately positive and negative,
and it follows that $n\alpha -m$ lies between the sum of the odd terms and the sum of the even terms.
To bound these terms,
we use
$d_\nu\leq a_{\nu +1}$.
Moreover,
$d_k\geq 1$,
hence $d_{k+1}\leq a_{k+2}-1$.
It follows that $n\alpha -m$ lies strictly between
$$
a_{k+1}(\alpha q_k-p_k)+a_{k+3}(\alpha q_{k+2}-p_{k+2})+a_{k+5}(\alpha q_{k+4}-p_{k+4})+\dots
$$
and
$$
(\alpha q_k-p_k)+(a_{k+2}-1)(\alpha q_{k+1}-p_{k+1})+a_{k+4}(\alpha q_{k+3}-p_{k+3})+\dots .
$$
Now
$
a_{\nu +1}(\alpha q_\nu -p_\nu )=\alpha (q_{\nu +1}-q_{\nu -1})-(p_{\nu +1}-p_{\nu -1}).
$
So both sums are telescopic.
The first sum immediately evaluates to $-\alpha q_{k-1}+p_{k-1}=(-1)^k/q'_k$.
The second sum equals $(\alpha q_{k}-p_{k})-(\alpha q_{k+1}-p_{k+1})-(\alpha q_{k}-p_{k})=(-1)^k/q'_{k+2}$.
\end{proof}

\subsection{Two Generators: the Bernoulli Flow}

Assume that $N=2$,
and let $f$ be defined as in~(\ref{E: f}) with weights $w_1$ and $w_2=\alpha w_1$,
for some irrational number $\alpha >1$.
We want to study the complex solutions to the equation $f(\omega )=0$ that lie close to the line $\Re s=D$.
First of all,
such solutions must have $e^{-w_1\omega}$ close to $e^{-w_1D}$,
so we take $\omega$ to be close to $D+2\pi iq/w_1$,
for an integer $q$.
Then we write $\alpha q=p+\err /(2\pi i)$,
for an integer $p$,
which we will specify below,
and $\omega =D+2\pi iq/w_1+\Domega$.
With these substitutions,
the equation $f(\omega )=0$ transforms to
$1 - e^{-w_1D}e^{-w_1\Domega } - e^{-w_2D}e^{-\err }e^{-w_2\Domega} = 0$.
This equation defines~$\Domega $ as a function of $\err $.

\begin{lemma}\label{L: approx}
Let\/ $w_1,w_2>0$ and\/ $\alpha =w_2/w_1>1;$
let\/ $D$ be such that\/ $e^{-w_1D}+ e^{-w_2D}=1,$
and let\/ $\Domega =\Domega (\err )$ be the function of\/ $\err ,$
defined implicitly by
\begin{gather}
e^{-w_1D}e^{-w_1\Domega }+ e^{-w_2D}e^{-\err }e^{-w_2\Domega }=1,
\end{gather}
and $\Domega (0)=0.$
Then $\Domega $ is analytic in~$\err ,$
in a disc of radius at least\/ $\pi$ around\/ $\err =0,$
with power series
\begin{align*}
\Domega (\err )&=-\frac{e^{-w_2D}}{f'(D)}\, x
+\frac{w_1^2e^{-w_1D} e^{-w_2D}}{2{f'(D)}^3} \, x^2 +O(x^{3}),
\qquad\mbox{as }\err \to 0.
\end{align*}
The coefficients of this power series are real.
The coefficient of\/ $x$ is negative and that of\/ $x^2$ is positive.
\end{lemma}

\begin{proof}
Define $y=y(\err )$ by $e^{-w_1D}y+e^{-w_2D}e^{-\err }y^\alpha =1$.
Then $y(0)=1$ and $w_1\Domega =-\log y$.
Since $y$ does not vanish,
it follows that if $y(\err )$ is analytic in a disc centered at $\err =0$,
then $\Domega $ will be analytic in that same disc.
Moreover,
$y$ is real-valued and positive when $\err $ is real,
so $\Domega $ is real-valued as well when $\err $ is real.
Further,
$y(\err )$ is locally analytic in $\err $,
with derivative
$$
y'(\err )=\frac{e^{-w_1D}y^\alpha e^{-x}}{e^{-w_1D}+\alpha e^{-w_2D}y^{\alpha -1}e^{-x}}.
$$
Hence there is a singularity at those values of $\err $ at which the denominator vanishes,
which is at $y=(\alpha /(\alpha -1))e^{w_1D}$ and
$
e^{-\err }=-\alpha^{-\alpha}(\alpha -1)^{\alpha -1}.
$
Since this value is negative,
the disc of convergence of the power series for $y(\err )$ is
$$
|x|<\left| -\alpha\log\alpha +(\alpha -1)\log (\alpha -1)+\pi i\right| ,
$$
which is a disc of radius at least $\pi$.
The first two terms of the power series for $\Domega (\err )$ are now readily computed.
\end{proof}

Applying this,
we find
\begin{gather}\label{E: omega}
\omega =D+2\pi i\frac{q}{w_1}-\frac{e^{-w_2D}}{f'(D)}\, x
+\frac{w_1^2e^{-w_1D} e^{-w_2D}}{2{f'(D)}^3} \, x^2
+O(\err^3),
\end{gather}
as $\err =2\pi i(q\alpha -p)\to 0$.
We view this formula as expressing $\omega$ as an initial approximation $D+2\pi iq/w_1$,
which is corrected by each term in the power series.
The first corrective term is in the imaginary direction,
as are all the odd ones,
and the second corrective term,
along with all the even ones,
are in the real direction.
The second term decreases the real part of $\omega$.

\begin{theorem}\label{T: cxd at q}
Let\/ $\alpha$ be irrational and let\/ $p_\nu$ and\/ $q_\nu$ be defined by~{\rm (\ref{E: p/q}).}
Let\/~$q$ be a positive integer,
and let\/ $q=\sum_{\nu =k}^\infty d_\nu q_\nu$ be the $\alpha$-adic expansion of\/ $q,$
as in Lemma\/~{\ref{L: adic expansion}.}
Assume $k\geq 2$ or $k=1$ and\/ $a_1\geq 2,$
and put\/ $p=\sum_{\nu =k}^\infty d_\nu p_\nu .$
Then there exists a complex dimension of\/ $\Flow_\weight$ at
\begin{equation}
\begin{aligned}
\omega &=D+2\pi i\frac{q}{w_1}-2\pi i\frac{e^{-w_2D}}{f'(D)}(q\alpha -p)\\
&\hspace{5mm}
-2\pi^2\frac{w_1^2e^{-w_1D} e^{-w_2D}}{{f'(D)}^3}(q\alpha -p)^2
+O\left( (q\alpha -p)^3\right) .
\end{aligned}
\end{equation}
The imaginary part of this complex dimension is approximately $2\pi iq/w_1,$
and its distance to the line $\Re s=D$ is at least\/ $C/q_{k+2}^{\prime 2},$
where $C=2\pi^2w_1^2e^{-(w_1+w_2)D}/{f'(D)}^3$ depends only on $w_1$ and\/ $w_2.$

Moreover,
$|\zeta_\weight (s)|\ll q^{\prime 2}_{k+2}$ around $s=D+2\pi iq/w_1$ on the line $\Re s=D,$
and\/ $|\zeta_\weight (s)|$ reaches a maximum of size $C'(q\alpha -p)^{-2},$
where $C'$ depends only on the weights $w_1$ and\/ $w_2.$
\end{theorem}

\begin{proof}
By Lemma~\ref{L: adic expansion},
the quantity $q\alpha -p$ lies between $(-1)^k/q'_{k+2}$ and $(-1)^k/q'_{k}$.
Under the given conditions on $k$,
$q'_k>q_k\geq 2$,
hence $\err =2\pi i(q\alpha -p)$ is less than $\pi$ in absolute value.
Then (\ref{E: omega}) gives the value of $\omega$.

Since the derivative of $f$ is bounded on the line $\Re s=D$,
this also implies that $f(s)$ reaches a minimum of order $(q\alpha -p)^{2}$ on an interval around
$s=D+2\pi iq/w_1$ on the line $\Re s=D$.
It follows that $|\zeta_\weight (s)|\ll q^{\prime 2}_{k+2}$ on the line $\Re s=D$,
with a maximum of order $(q\alpha -p)^{-2}$.
\end{proof}

We obtain more precise information when $q=q_k$.

\begin{theorem}\label{T: cxd at qk}
For every\/ $k\geq 0$ {\rm (}or $k\geq 1$ if\/ $a_1=1),$
there exists a complex dimension $\omega$ of\/ $\Flow_\weight$ of the form
\begin{equation}
\begin{aligned}
\omega &=D+2\pi i\frac{q_k}{w_1}-2\pi i(-1)^k\frac{e^{-w_2D}}{f'(D)q'_{k+1}}
-2\pi^2w_1^2\,\frac{e^{-(w_1+w_2)D}}{{f'(D)}^3q_{k+1}^{\prime 2}}
+O\left( q_{k+1}^{\prime -3}\right) ,
\end{aligned}
\end{equation}
as $k\to\infty .$

Moreover,
$|\zeta_\weight (s)|\ll q^{\prime 2}_{k+1}$ around $s=D+2\pi iq_k/w_1$ on the line $\Re s=D,$
and\/ $|\zeta_\weight (s)|$ reaches a maximum of size $C'q_{k+1}^{\prime 2},$
where $C'$ is as in Theorem~\ref{T: cxd at q}.
\end{theorem}

\begin{proof}
Put $p=p_k$.
Then $\err =2\pi i(-1)^k/q'_{k+1}$,
which is less than $\pi$ in absolute value.
The rest of the proof is the same as in the proof of Theorem~\ref{T: cxd at q}.
\end{proof}

\begin{definition}
A domain in the complex plane containing the line $\Re s=D$ is a {\em dimension free region\/} for
the flow $\Flow_\weight$ if the only pole of $-\zeta_\weight'/\zeta_\weight$ in that region is $s=D$.
\end{definition}

\begin{corollary}
Assume that the coefficients $a_0,a_1,\dots$ of\/ $\alpha$ are bounded by $M.$
Put $B=\pi^4e^{-(w_1+w_2)D}/(2{f'(D)}^3).$
Then $\Flow_\weight$ has a dimension free region of the form
\begin{gather}
\left\{\sigma +it\in\Com\colon\sigma >D-\frac{B}{M^{2}t^{2}}\right\} .
\end{gather}
The function $-\zeta_\weight'/\zeta_\weight$ satisfies hypotheses {\rm ({\bf H$_1$})} and\/
{\rm (\bf H$_2$\rm )} with\/ $\kappa =2.$

More generally,
let\/ $a\colon\Real^+\rightarrow [1,\infty )$ be a function such that the coefficients $\{ a_k\}_{k=0}^\infty$
of the continued fraction of\/ $\alpha$ satisfy $a_{k+1}\leq a(q_k)$ for every $k\geq 0.$
Then $\Flow_\weight$ has a dimension free region of the form
\begin{gather}
\left\{\sigma +it\in\Com\colon\sigma >D-\frac{B}{t^{2}a^{2}(tw_1/(2\pi ))}\right\} .
\end{gather}
If\/ $a$ grows at most polynomially,
then $-\zeta_\weight'/\zeta_\weight$ satisfies hypotheses {\rm ({\bf H$_1$})} and\/
{\rm (\bf H$_2$\rm )} with\/ $\kappa$ such that\/ $t^\kappa\geq t^2a^2(tw_1/(2\pi )).$
\end{corollary}

\begin{proof}
This follows from Theorem~\ref{T: cxd at qk},
if we note that for $t=2\pi q_k/w_1$,
we have $q'_{k+1}=\alpha_{k+1}q'_k\leq 2a(q_k)q'_k\leq 4a(q_k)q_k$.
So the complex dimension close to $D+it$ is located at
$D+i(t+O(q_{k+1}^{\prime -1}))-(w_1^2/\pi^2)Bq_{k+1}^{\prime -2}+O(q_{k+1}^{\prime -4})$,
where the orders denote real-valued functions.
The real part of this complex dimension is less than $D-Bt^{-2}a^{-2}(tw_1/(2\pi ))$.
\end{proof}

This has the following consequence for the Prime Orbit Theorem:

\begin{theorem}[Prime Orbit Theorem with Error Term, for Bernoulli flows]\label{T: nl PNT}
$\ $
Let\/ $\alpha =w_2/w_1$ have bounded coefficients in its continued fraction.
Then 
\begin{gather}
\psi_\weight (x)=\frac{x^D}{D}+O\left( x^D\left(\frac{\log\log x}{\log x}\right)^{1/4}\right) ,
\end{gather}
as $x\to\infty .$

If\/ $\alpha$ is `polynomially approximable',
with coefficients in its continued fraction satisfying $a_{k+1}\leq a(q_k),$
for some increasing function $a$ with\/ $a(x)=O(x^l),$
as $x\to\infty ,$
then
\begin{gather}
\psi_\weight (x)=\frac{x^D}{D}
+O\left( x^D \left(\frac{\log\log x}{\log x}\right)^{\frac{1}{4l+4}}\right) ,
\end{gather}
as $x\to\infty .$
\end{theorem}

The proof will be given in the next section,
see Theorem~\ref{T: nl PNT N}.

\begin{example}\label{E: cxd of GF}
For the Golden flow,
we have $\alpha =\phi$ and $w_1=\log 2$ (see Example~\ref{E: golden flow}).
The continued fraction of $\phi$ is $[1,1,1,\dots ]$,
hence $q_k=F_{k+1}$,
the $(k+1)$-st Fibonacci number,
and $q'_k=\phi^k$.
Numerically,
we find $D\approx .7792119034$ and the following approximation of $\Domega (\err )$:
\begin{align*}
&-.47862\, x+.08812\, x^2+.00450\, x^3-.00205\, x^4-.00039\, x^5+.00004\, x^6+\dots .
\end{align*}
For every $k$,
we find a complex dimension close to $D+2\pi iq_k/\log 2$.
For example,
$q_9=55$,
and we find a complex dimension at $D-.00023 + 498.58 i$.
More generally,
for numbers like $q=55+5$ or $q=55-5=34+13+3$,
we find a complex dimension close to $D+2\pi iq/\log 2$,
in this case respectively at $D-.023561 + 543.63 i $ and at $D-.033919 + 453.53 i $.
In both these cases,
the distance to the line $\Re s=D$ is comparable to the distance of the complex dimension for $q=5$ to this
line,
which is located at $D-.028499 + 45.05 i$.
See Figure~\ref{F: golden global},
where the markers are at the Fibonacci numbers.
The pattern persists for other complex dimensions as well.
Indeed,
every complex dimension repeats itself according to the Fibonacci numbers.
\end{example}

\subsection{More than Two Generators}

The following lemma replaces the continued fraction construction.

\begin{lemma}\label{L: multi approx}
Let\/ $w_1,w_2,\dots ,w_N$ be weights such that at least one ratio $w_j/w_k$ is irrational.
Then for every $Q>1,$
there exist integers $q<Q^{N-1}$ and\/ $p_j$ such that $|qw_j-p_jw_1|\leq w_1/Q$ for $j=1,\dots ,N.$
In particular,
$|qw_j-p_jw_1|<w_1q^{-1/(N-1)}$ for $j=1,\dots ,N.$
\end{lemma}

\begin{remark}
Note that the condition implies that at least one ratio $w_j/w_1$ is irrational.
Also,
$|qw_j-p_jw_1|\neq 0$ when $w_j/w_1$ is irrational,
so $q\to\infty$ when $Q\to\infty$.

The construction of such integers $q$ and $p_j$ is much less explicit than for $N=2$,
since there does not exist a continued fraction algorithm for simultaneous approximation.\footnote
{However,
the L$^3$-algorithm can be used as a substitute for the continued fraction algorithm.
We thank H.~W.~Lenstra, Jr.\ for guiding us to the following information:
The L$^3$-algorithm~\cite{LLL} can be used to find fractions $p_j/q$ that approximate $w_j/w_1$ for
$j=1,\dots ,N$.
This algorithm works in polynomial time,
like the continued fraction algorithm,
but it does not give the best possible value for $q$ (given a certain error of approximation).
The problem of finding the best value for $q$ is NP-complete~\cite{jL82}.
(See also~\cite{GLS}.)}
The number $Q$ plays the role of $q'_{k+1}$ in Theorem~\ref{T: cxd at qk} above.
In particular,
if~$q$ is often much smaller than $Q$,
then $w_1,\dots ,w_N$ is well approximable by rationals,
and we find a small dimension free region.
\end{remark}

Again,
we are looking for a solution of $f(\omega )=0$ close to $s=D+2\pi iq/w_1$,
where~$f$ is defined by~(\ref{E: f}).
We write $\omega =D+2\pi iq/w_1+\Domega$ and $w_jq=w_1p_j+w_1x_j/(2\pi i)$.
For $j=1$,
we take $p_1=q$ and consequently $x_1=0$.
In general,
$x_j=2\pi i(qw_j/w_1-p_j)$.
Then $f(\omega )=0$ is equivalent to $1-\sum_{j=1}^Ne^{-w_jD}e^{-x_j-w_j\Domega}=0$.

The following lemma is the several variable analogue of Lemma~\ref{L: approx}.
However,
in this case we do not know the radius of convergence with respect to each of the variables involved.

\begin{lemma}
Let\/ $w_1\leq w_2\leq\dots\leq w_N,$
let\/ $D$ be such that\/ $\sum_{j=1}^Ne^{-w_jD}=1,$
and let\/ $\Domega =\Domega (x_2,\dots ,x_N)$ be implicitly defined by
\begin{gather}
\sum_{j=1}^Ne^{-w_jD}e^{-x_j-w_j\Domega}=1,
\end{gather}
and\/ $x_1=0.$
Then $\Domega$ is analytic in $x_2,\dots ,x_N,$
with power series
\begin{equation}
\begin{aligned}
\Domega =
&-\sum_{j=2}^N\frac{e^{-w_jD}}{f'(D)}\, x_j
+\frac{1}{2}\sum_{j=2}^N\frac{e^{-w_jD}}{f'(D)}\, x_j^2\\
&-\frac{1}{2}\sum_{j,k=2}^N\left(\frac{f''(D)}{{f'(D)}^3}+\frac{w_j+w_k}{{f'(D)}^2}\right)
e^{-(w_j+w_k)D}x_jx_k
+O\biggl(\sum_{j=2}^N|x_j|^3\biggr).
\end{aligned}
\end{equation}
This power series has real coefficients.
The terms of degree two form a positive definite quadratic form.
\end{lemma}

\begin{proof}
The positive definiteness follows from the fact that the complex dimensions lie to the left of
$\Re s=D$,
see Theorem~\ref{T: struc cxd}.
It can also be verified directly.
\end{proof}

Applying this,
we find
\begin{equation}\label{E: omegaN}
\begin{aligned}
\omega &=D+2\pi i\frac{q}{w_1}-\sum_{j=2}^N\frac{e^{-w_jD}}{f'(D)}\, x_j
+\frac{1}{2}\sum_{j=2}^N\frac{e^{-w_jD}}{f'(D)}\, x_j^2\\
&\hspace{5mm}
-\frac{1}{2}\sum_{j,k=2}^N\left(\frac{f''(D)}{{f'(D)}^3}+\frac{w_j+w_k}{{f'(D)}^2}\right)
e^{-(w_j+w_k)D}x_jx_k
+O\biggl(\sum_{j=2}^N|x_j|^3\biggr) ,
\end{aligned}
\end{equation}
where $x_j=2\pi i(qw_j/w_1-p_j)$.
Again,
this formula expresses $\omega$ as an initial approximation $D+2\pi iq/w_1$,
which is corrected by each term in the power series.
The corrective terms of degree one are again in the imaginary direction,
as are all the odd degree ones,
and the corrective terms of degree two,
along with all the even ones,
are in the real direction.
The degree two terms decrease the real part of $\omega$.

\begin{theorem}
Let\/ $N\geq 2$ and let\/ $w_1,\dots ,w_N$ be weights.
Let $Q$ and $q$ be as in Lemma~\ref{L: multi approx}.
Then $\Flow_\weight$ has a complex dimension close to $D+2\pi iq/w_1$ at a distance of at most\/ $O(Q^{-2})$
from the line $\Re s=D,$
as $Q\to\infty .$
The function $|\zeta_\weight'/\zeta_\weight |$ reaches a maximum of order $Q^2.$
\end{theorem}

\begin{proof}
Again,
the numbers $x_j$ are purely imaginary,
so the corrective terms of degree~$1$ (and of every odd degree) give a correction in the imaginary direction,
and only the corrective terms of even degree will give a correction in the real direction.
Since $|x_j|<2\pi /Q$,
the theorem follows.
\end{proof}

\begin{corollary}
The best dimension free region that $\Flow_\weight$ can have is of size
\begin{gather}
\left\{ \sigma +it\colon\sigma\geq D-O\left( t^{-2/(N-1)}\right)\right\} .
\end{gather}
The implied constant depends only on $w_1,\dots ,w_N.$

If\/ $w_1,\dots ,w_N$ is `$a$-approximable',
then the dimension free region has the form
\begin{gather}
\left\{ \sigma +it\colon\sigma\geq D-O\left( a^{-2}(w_1t/(2\pi ))t^{-2/(N-1)}\right)\right\} ,
\end{gather}
where $a\colon [1,\infty )\rightarrow\Real^+$ is an increasing function such that for every integer $q\geq 1,$
$|qw_j-p_jw_1|\geq (w_1/a(q))q^{-1/(N-1)}$ for $j=1,\dots ,N.$
\end{corollary}

This has the following consequence for the Prime Orbit Theorem:

\begin{theorem}[Prime Orbit Theorem with Error Term]\label{T: nl PNT N}
Suppose\/ $w_1,\dots ,w_N$ are badly approximable,
in the sense that\/ $|qw_j-p_jw_1|\gg q^{-1/(N-1)}$ for $j=1,\dots ,N$ and every $q\geq 1.$
Then
\begin{gather}
\psi_\weight (x)=\frac{x^D}{D}+O\left( x^D\left(\frac{\log\log x}{\log x}\right)^{\frac{N-1}4}\right) ,
\end{gather}
as $x\to\infty .$

If\/ $w_1,\dots ,w_N$ is `polynomially approximable',
in the sense that\/ $|qw_j-p_jw_1|\geq (w_1/a(q))q^{-1/(N-1)}$ for $j=1,\dots ,N$ and every $q\geq 1,$
for some increasing function $a$ on $[1,\infty )$ such that\/ $a(x)=O(x^l)$ as $x\to\infty ,$
then
\begin{gather}
\psi_\weight (x)
=\frac{x^D}{D}+O\left( x^D \left(\frac{\log\log x}{\log x}\right)^{\frac{N-1}{4l(N-1)+4}}\right) ,
\end{gather}
as $x\to\infty .$
\end{theorem}

\begin{proof}[Proof of Theorems~\ref{T: nl PNT} and~\ref{T: nl PNT N}]\label{P: nl PNT N}
We apply the pointwise explicit formula at level $k=2$ (see Theorem~\ref{T: pointwise}) to obtain
$$
\psi_\weight^{[2]}(x)=\frac{x^{D+1}}{D(D+1)}
+\sum_{\omega\in\Dimensions_\weight\backslash\{ D\}}\frac{x^{\omega +1}}{\omega (\omega +1)}+R^{[2]} (x).
$$
The error term is estimated by $R^{[2]}(x)=O(x^{D+1-c})$ for some positive $c$.
We will estimate the sum by an argument which is classical in the theory of the Riemann zeta function and the
Prime Number Theorem,
under the assumptions that the $\omega$ have a linear density,
and that every $\omega =\sigma +it$ satisfies $\sigma\leq D-Ct^{-\rho}$ for some positive number $\rho$.
We then obtain Theorem~\ref{T: nl PNT N} by taking $\rho =2/(N-1)+2l$,
and Theorem~\ref{T: nl PNT} corresponds to the case when $N=2$.

The sum $\sum_\omega \frac{x^{\omega +1}}{\omega (\omega +1)}$ is absolutely convergent.
We split this sum into the parts with $|\Im\omega |>T$ and with $|\Im\omega |\leq T$.
Put $A=\sum_\omega\left|{\omega (\omega +1)}\right|^{-1}$.
From the fact that the complex dimensions have a linear density,
it follows that there exists a constant $B$ is such that
$
\sum_{|\Im\omega |\geq T}\left|{\omega (\omega +1)}\right|^{-1}\leq{B}/{T}
$
for every $T$.
Then $\left|\sum_\omega \frac{x^{\omega +1}}{\omega (\omega +1)}\right|\leq Ax^{D+1-CT^{-\rho }}+Bx^{D+1}/T$.
For $T=(\rho C\log x/\log\log x)^{1/\rho }$,
we find
$$
\left|\sum_\omega \frac{x^{\omega +1}}{\omega (\omega +1)}\right|
=O\left( x^{D+1}\left(\frac{\log\log x}{\log x}\right)^{1/\rho}\right) .
$$

We then apply a Tauberian argument to deduce a similar error estimate for $\psi_\weight (x)$;
see~\cite[p.~64]{Ingham}.
Let $h=x(\log\log x/\log x)^{1/(2\rho )}$.
Thus
$$
\psi_\weight (x)\leq\frac{1}{h}\int_x^{x+h}\psi_\weight (t)\, dt
=\frac{\psi_\weight^{[2]}(x+h)-\psi_\weight^{[2]}(x)}{h}.
$$
Now $\frac{(x+h)^{D+1}-x^{D+1}}{hD(D+1)}=x^D/D+O(x^{D-1}h)=x^D/D+x^DO((\log\log x/\log x)^{1/(2\rho )})$.
Further,
$O(x^{D+1}(\log\log x/\log x)^{1/\rho }/h)=x^DO((\log\log x/\log x)^{1/(2\rho )})$.
\end{proof}

\begin{remark}
Note that by using the Tauberian argument,
we lose a factor two in the exponent.
Indeed,
the estimate
$$
\sum_{\omega\in\Dimensions_\weight\backslash\{ D\}}\frac{x^\omega}{\omega}+R (x)
=O\left( x^D \left(\frac{\log\log x}{\log x}\right)^{\frac{N-1}{2l(N-1)+2}}\right)
$$
holds distributionally.
\end{remark}

\begin{remark}
If $a(q)$ grows more than polynomially,
we obtain a bound of the form $x^D/a^{\rm inv}(\log x)$ for the error in the Prime Orbit Theorem,
where $a^{\rm inv}$ is the inverse function of $a$.
\end{remark}

\subsection{Conclusion}

If $l>0$ in the exponent $(N-1)/(4l(N-1)+4)$ of $\log x$ in the error term of Theorems~\ref{T: nl PNT}
and~\ref{T: nl PNT N},
then the error term is independent of~$N$,
essentially of order $x^D(\log x)^{-1/4l}$ (ignoring the factor of $\log\log x$).
Thus,
if the weights are well approximable,
the error term is never better than $x^D$ divided by a fixed power of the logarithm of $x$.
On the other hand,
when $l=0$,
that is,
roughly speaking,
when the weights are never close to rational numbers,
then the error term is essentially of order $x^D(\log x)^{-(N-1)/4}$.
Hence,
the larger $N$,
the smaller the error term in that case.

We may compare this,
somewhat superficially in view of the Riemann Hypothesis,
with the situation of the Riemann zeta function.
In view of Example~\ref{E: Riemann-von Mangoldt},
the weights are $w_p=\tweight (p)=\log p$,
for each prime number $p$,
and there are infinitely many of them.
Since it is expected that $\{\log p\}_{p\text{:prime}}$ is badly approximable,
one expects an error term of order ``$x^D(\log x)^{-\infty}$''.
Indeed,
in~(\ref{E: classical PNT}),
we have $e^{-c\sqrt{\log x}}=O\left( (\log x)^{-N}\right)$ for every $N>0$.
The corresponding pole free region has width $A/\log t$ at height $t$ (see~\cite[Theorem~19]{Ingham}),
which is ``$t^{-1/\infty}$''.
This lends credibility to the conjecture that $\{\log p\}_{p\text{:prime}}$ is badly approximable by rational
numbers.
\medskip

{\bf Acknowledgement.}
The authors wish to thank Gabor Elek for helpful conversations about dynamical zeta
functions.

\bibliographystyle{amsalpha}

\end{document}